\documentclass{elsart}%
\usepackage{amsmath}
\usepackage{graphicx}
\usepackage{float}
\usepackage{indentfirst}
\usepackage{verbatim}
\usepackage{algorithm}
\usepackage{algpseudocode}
\usepackage{float}
\usepackage{multicol}
\usepackage{listings}
\usepackage{subfigure}
\usepackage{epsfig}
\usepackage{amsfonts}
\usepackage{amssymb}
\usepackage{epstopdf}
\usepackage{color}%
\begin{document}
\begin{frontmatter}
\title{Highly Accurate Nystr\"{o}m Volume Integral Equation Method for the Maxwell equations for 3-D Scatters}
\author[UNCC]{Duan Chen}
\author[UNCC]{Wei Cai}
\author[UNCC]{Brian Zinser}
\address[UNCC]{Department of Mathematics and Statistics, University of
North Carolina at Charlotte, Charlotte, NC 28223, USA}
\bigskip
{\bf Suggested Running Head:}
\\
Highly Accurate Nystr\"{o}m Volume Integral Equations for the Maxwell equations
\\
\bigskip
{\bf Corresponding Author: }
\\
Prof. Wei Cai \\
Department of Mathematics and Statistics, \\
University of North Carolina at
Charlotte, \\
Charlotte, NC 28223-0001 \\
Phone: 704-687-0628, Fax: 704-687-6415, \\
Email: wcai@uncc.edu
\newpage
\begin{abstract}
In this paper, we develop highly accurate Nystr\"{o}m methods for the volume integral equation (VIE) of the Maxwell equation for 3-D scatters. The method is based on a formulation of the VIE equation where the Cauchy principal value of the dyadic Green's function can be computed accurately for a finite size exclusion volume with some explicit corrective integrals of removable singularities. Then, an effective interpolated quadrature formula for tensor product Gauss quadrature nodes in a cube is proposed to handle the hyper-singularity of integrals of the dyadic Green's function. The proposed high order Nystr\"{o}m VIE method is shown to have high accuracy and demonstrates $p$-convergence for computing the electromagnetic scattering of cubes in $R^3$.
\end{abstract}
\end{frontmatter}

\bigskip

\section{Introduction}

Electromagnetic (EM) wave scattering in the presence of (random)
microstructure on material interfaces has a wide range of applications. For
example, the interaction of light between the surface plasmon, a collective
fluctuation of electron density under the background of the positive nucleus
charges, on the metallic surfaces produces a so-called surface plasmon
polariton (SPP) \cite{Atwater:2007,Raether:1988}, which has important
applications in solar cells \cite{Atwater:2010}, meta-materials, and
super-resolution imaging devices \cite{Pendry:2000,Fu:2010}. Additionally,
surface enhanced Raman scattering (SERS) \cite{Raman:1928} is also closely
related to the excitation of surface plasmons on rough or nano-pattern
surfaces by incident light and is an extremely useful tool in finger-printing
the chemical components of a molecule, single molecule detector, detection of
DNA studies, and bio-sensor, etc \cite{Hering:2008}. In all these
applications, it is critical to have highly accurate and efficient numerical
methods for computer simulations of the EM scattering in microstructures.

In this paper, we will present a high order Nystr\"{o}m volume integral
equation (VIE) method for the time harmonic Maxwell equations using dyadic
Green's functions $\overline{\mathbf{G}}_{\mathbf{E}}(\mathbf{r}^{\prime
},\mathbf{r})$. In most of the applications, the scatter is embedded into
either a homogeneous or layered media and a Lippmann-Schwinger type of VIE,
which is a well-conditioned second kind of Fredhohm integral equation, can be
derived for the regions occupied by the scatter while the dyadic Green's
function $\overline{\mathbf{G}}_{\mathbf{E}}(\mathbf{r}^{\prime},\mathbf{r})$
will ensure that the scattering field, expressed in terms of equivalent
current sources inside the scatter, satisfies interface conditions along the
layered interfaces as well as the Sommerfeld radiation conditions at $\infty.$
In the VIE formulation, the electric field inside the scatter will use the
Cauchy Principal Value (CPV or simply p.v.) associated with the dyadic Green's
function as indicated in (\ref{per-197}). Therefore, one of the most difficult
issue is how to compute \textit{accurately and efficiently }the CPV for the
dyadic Green's function which has an $O(\frac{1}{R^{3}})$ singularity, i.e.
\begin{equation}
{\small \text{ p.v.}\int_{\Omega}\mathrm{d}\mathbf{r}^{\prime}\text{
}\mathrm{i}\omega\Delta\epsilon(\mathbf{r}^{\prime})\mathbf{E}(\mathbf{r}%
^{\prime})\cdot\overline{\mathbf{G}}_{\mathbf{E}}(\mathbf{r}^{\prime
},\mathbf{r})=\lim_{\delta\rightarrow0}\int_{\Omega\backslash V_{\delta}%
}\mathrm{d}\mathbf{r}^{\prime}\text{ }\mathrm{i}\omega\Delta\epsilon
(\mathbf{r}^{\prime})\mathbf{E}(\mathbf{r}^{\prime})\cdot\overline{\mathbf{G}%
}_{\mathbf{E}}(\mathbf{r}^{\prime},\mathbf{r}),} \label{cpv}%
\end{equation}
where ${\small \Delta\epsilon(\mathbf{r}^{\prime})}$ describes the scatter
dielectric constant deviating from the surrounding background material. As the
CPV is defined through a limiting process of diminishing size $\delta$ of the
exclusion volume, in practical computation, a finite $\delta$ has to be taken,
namely, a small, finite and fixed $\delta>0$ for the exclusion $V_{\delta}$ is
selected, then the following approximation is taken
\begin{equation}
\text{ p.v.}\int_{\Omega}\mathrm{d}\mathbf{r}^{\prime}\text{ }\mathrm{i}%
\omega\Delta\epsilon(\mathbf{r}^{\prime})\mathbf{E}(\mathbf{r}^{\prime}%
)\cdot\overline{\mathbf{G}}_{\mathbf{E}}(\mathbf{r}^{\prime},\mathbf{r}%
)\approx\int_{\Omega\backslash V_{\delta}}\mathrm{d}\mathbf{r}^{\prime}\text{
}\mathrm{i}\omega\Delta\epsilon(\mathbf{r}^{\prime})\mathbf{E}(\mathbf{r}%
^{\prime})\cdot\overline{\mathbf{G}}_{\mathbf{E}}(\mathbf{r}^{\prime
},\mathbf{r}). \label{eqn:numericalpv}%
\end{equation}

Therefore, we are faced with two issues in the implementation of the VIE:

Firstly, what size of $\delta$ should be taken? What is the effect of error by
using a finite $\delta$ in the calculation of the CPV? For any finite
$\delta,$ there will be a truncation error which we will call correction
terms
\begin{align}
\text{ p.v.}\int_{\Omega}\mathrm{d}\mathbf{r}^{\prime}\text{ }\mathrm{i}%
\omega\Delta\epsilon(\mathbf{r}^{\prime})\mathbf{E}(\mathbf{r}^{\prime}%
)\cdot\overline{\mathbf{G}}_{\mathbf{E}}(\mathbf{r}^{\prime},\mathbf{r})  &
=\int_{\Omega\backslash V_{\delta}}\mathrm{d}\mathbf{r}^{\prime}\text{
}\mathrm{i}\omega\Delta\epsilon(\mathbf{r}^{\prime})\mathbf{E}(\mathbf{r}%
^{\prime})\cdot\overline{\mathbf{G}}_{\mathbf{E}}(\mathbf{r}^{\prime
},\mathbf{r})\nonumber\\
&  +\text{ correction terms}. \label{correction}%
\end{align}
The existence of these correction terms and their magnitude will limit the
accuacy of the VIE solution if they are not explicitly included in the
numerical solution process. The correction terms were derived by Fikioris
\cite{Fik:1965} using the mixed potential formulation of the electric field.
In this paper, we will re-derive the VIE similar to those in \cite{Fik:1965},
however, in a more succinct manner and the final formula is more suitable for
numerical implementations. Previous work on how to handle the singular
integrals for VIE method include singularity substraction \cite{Kottmann:2000}%
, locally corrected Nystr\"{o}m scheme \cite{Liu:2001}, direct integration of
the singularity \cite{Tong:2010},\ etc.

Secondly, for the selected $\delta$, how to find an accurate integration
quadrature formula to compute the resulting integration over the domain
$\Omega\backslash V_{\delta}$ with a highly singular integrand?

It is the objective of this paper to address these two issues and find easily
implementable solutions. To address the first issue, we will derive the VIE
equation using vector and scale potentials such that the Cauchy principal
value can be computed as in (\ref{correction}) with explicit expression for
the correction terms. On the other hand, to address the second issue, we will
construct special quadrature weights for the tensor product Gauss quadrature
nodes in a reference element $\Omega$ (taken to be a cube in this paper) by
using an interpolation approach. In this approach, first a brute force
computation of the integral, using Gauss quadrature in polar coordinate
centered at the singularity, will be done to a given accuracy, albeit
involving large number of values of the integrand. Then, realizing the
integrand in the VIE matrix entries, except for the singular denominators
involving $R^{k},1\leq k\leq3$, are smooth functions, which can be in fact
accurately interpolated using values on the tensor product Gauss nodes inside
the cubic domain $\Omega.$ Then, the brute-force integration formula can be
converted into a new integration weights for the tensor product Gauss nodes.
The new integration formula can be tabulated for integrating general
functions. Finally, the Nystr\"{o}m collocation method is used to discretize
the VIE with a high order accuracy.

The rest of the paper is organized as follows: Section 2 presents the
formulation of a VIE where the CPV can be computed with a finite exclusion
volume accurately. Then, numerical algorithms, including the Nystr\"{o}m
collocation method, efficient quadrature formula and numerical implementation
are given in Section \ref{sec:methods}. Section \ref{sec:results} includes
various numerical simulation tests such as the accuracy of Cauchy principal
value computation, $\delta$-independence of matrix entries and $p$-refinement
convergence for the VIE. The paper ends with a conclusion in Section
\ref{sec:conclusion}.

\bigskip\ 

\section{Volume Integral equations for the Maxwell equations}

%\bigskip

\subsection{VIE and Cauchy principal values}

In this section, we will follow \cite{wcai:2013} to show briefly how the VIE
for the Maxwell equations can be derived using a vector form second Green's
identity for the following vector wave equations,
\begin{equation}
\mathcal{L}\mathbf{E}(\mathbf{r})\text{ }\mathbf{-}\text{ }\omega^{2}%
\epsilon(\mathbf{r})\mathbf{E}(\mathbf{r})=-\mathrm{i}\omega\mathbf{J}%
_{e}(\mathbf{r}),\text{ \ \ }\mathbf{r}\in\mathbb{R}^{3}\backslash(\Sigma
\cup\partial\Omega), \label{per-171}%
\end{equation}
where $\omega$ is the frequency, $\mu$ is the permeability, and  $\Sigma$ consists of any possible interfaces of the background medium in
case of a layered material,
\[
\mathcal{L}=\nabla\times\frac{1}{\mu}\nabla\times,
\]
and $\mathbf{J}_{e}(\mathbf{r})$ is the far-field source (assumed to be away
from the layered structure), which produces the incident waves impinging on
the layered structure from the top, i.e.,%
\begin{equation}
\mathbf{E}^{\text{inc}}(\mathbf{r})=-\mathrm{i}\omega\mu(\mathbf{r}%
)\int_{\mathbb{R}^{3}}\overline{\mathbf{G}}_{E}(\mathbf{r},\mathbf{r}^{\prime
})\cdot\mathbf{J}_{e}(\mathbf{r}^{\prime})\mathrm{d}\mathbf{r}^{\prime},
\label{per-173}%
\end{equation}
and $\overline{\mathbf{G}}_{E}(\mathbf{r},\mathbf{r}^{\prime})$ is the dyadic
Green's function for the layered media. A scatter $\Omega$ is characterized by
a different dielectric constant from the layered background dielectrics, i.e.,%
\begin{equation}
\epsilon(\mathbf{r})=\epsilon_{L}(\mathbf{r})+\Delta\epsilon(\mathbf{r}),
\label{per-175}%
\end{equation}
where $\Delta\epsilon(\mathbf{r})=0,\mathbf{r\notin}$ $\Omega.$ Then,
(\ref{per-171}) can be rewritten as%
\begin{equation}
\mathcal{L}\mathbf{E}(\mathbf{r})\text{ }\mathbf{-}\text{ }\omega^{2}%
\epsilon_{L}(\mathbf{r})\mathbf{E}(\mathbf{r})=-\mathrm{i}\omega
\mathbf{J}(\mathbf{r}), \label{per-177}%
\end{equation}
where%
\begin{equation}
\mathbf{J}(\mathbf{r})=\mathbf{J}_{e}(\mathbf{r})+\mathbf{J}_{\text{eq}%
}(\mathbf{r}), \label{per-178}%
\end{equation}
and the equivalent current source $\mathbf{J}_{\text{eq}}(\mathbf{r})$ is
defined to reflect the existence of the scatter $\Omega$:%
\begin{equation}
\mathbf{J}_{\text{eq}}(\mathbf{r})=\mathrm{i}\omega\Delta\epsilon
(\mathbf{r})\mathbf{E}(\mathbf{r}). \label{per-179}%
\end{equation}

Let us consider any interior point inside the scatter, i.e., $\mathbf{r}%
^{\prime}\in\Omega$ and a small volume $V_{\delta}=V_{\delta}(\mathbf{r}%
^{\prime})\subset\Omega$ centered at $\mathbf{r}^{\prime}.$ The dyadic Green's
function $\overline{\mathbf{G}}_{E}(\mathbf{r},\mathbf{r}^{\prime})$ is
defined by
\begin{equation}
\mathcal{L}\overline{\mathbf{G}}_{E}(\mathbf{r},\mathbf{r}^{\prime})\text{{}%
}\mathbf{-}\text{ }\omega^{2}\epsilon_{L}(\mathbf{r})\overline{\mathbf{G}}%
_{E}(\mathbf{r},\mathbf{r}^{\prime})=\frac{1}{\mu(\mathbf{r})}\overline
{\mathbf{I}}\delta(\mathbf{r}-\mathbf{r}^{\prime}),\text{ }\mathbf{r}%
\in\mathbb{R}^{3}.\label{per-185}%
\end{equation}
In the case of a homogeneous medium, we have%
\begin{equation}
\overline{\mathbf{G}}_{\mathbf{E}}\mathbf{(r},\mathbf{r^{\prime}%
)=\overline{\mathbf{G}}_{E}(r^{\prime}},\mathbf{r)=}\left(  \overline
{\mathbf{I}}\text{ }\mathbf{+}\text{ }\frac{1}{k^{2}}\nabla\nabla\right)
g(\mathbf{r},\mathbf{r^{\prime}}),\label{eqn3-34}%
\end{equation}
where $k^{2}=\omega^{2}\epsilon_{i}\mu,$ $\epsilon_{i}=\epsilon_{L}%
(\mathbf{r}),\mu$ is the permeability of the medium, and
\begin{equation}
g(\mathbf{r},\mathbf{r^{\prime}})=\frac{1}{4\pi}\frac{\mathrm{e}%
^{-\mathrm{i}kR}}{R},\qquad R=|\mathbf{r}-\mathbf{r}^{\prime}|.\label{eqn3-35}%
\end{equation}

Next, on multiplying (\ref{per-177}) by $\overline{\mathbf{G}}_{E}%
(\mathbf{r},\mathbf{r}^{\prime})$ and (\ref{per-185}) by $\mathbf{E}%
(\mathbf{r})$ and forming the difference, and then integrating over the domain
$\mathbb{R}^{3}\backslash V_{\delta},$ after some manipulation
\cite{wcai:2013}, we arrive at the following equation (after switching
$\mathbf{r}$ and $\mathbf{r}^{\prime}$):%
\begin{align}
&  -\mathrm{i}\omega\mu(\mathbf{r})\int_{\mathbb{R}^{3}\backslash V_{\delta}%
}\mathrm{d}\mathbf{r}^{\prime}\text{ }\overline{\mathbf{G}}_{E}(\mathbf{r}%
,\mathbf{r}^{\prime})\cdot\mathbf{J}(\mathbf{r}^{\prime})-\mu(\mathbf{r}%
)\int_{{S}_{\delta}}\ \mathrm{d}s^{\prime}\left[  \mathrm{i}\omega\text{{}%
}\overline{\mathbf{G}}_{E}(\mathbf{r},\mathbf{r}^{\prime})\cdot\left(
\mathbf{n}\times\mathbf{H}(\mathbf{r}^{\prime})\right)  \right. \nonumber\\
&  {-}\text{ }\left.  \frac{1}{\mu(\mathbf{r}^{\prime})}\nabla\times
\overline{\mathbf{G}}_{E}(\mathbf{r},\mathbf{r}^{\prime})\cdot\left(
\mathbf{n\times E}(\mathbf{r}^{\prime})\right)  \right]  =\mathbf{0}%
,\qquad\qquad\mathbf{r}\in\Omega, \label{per-186}%
\end{align}
where $S_{\delta}=\partial V_{\delta}(\mathbf{r}),$ $\mathbf{n}$ is the normal
of $S_{\delta}$ pointing out of $V_{\delta}(\mathbf{r})$.

As $\delta\rightarrow0$, the first integral will approach the Cauchy principal
value of a singular integral, while the surface integrals will in fact depend
on the geometric shape of the volume $V_{\delta}.$

In order to estimate the surface integrals, we have the following asymptotics
for small $kR\ll1$:
\begin{align}
\overline{\mathbf{G}}_{\mathbf{E}}\mathbf{(r},\mathbf{r^{\prime})}  &
=\frac{1}{4\pi k^{2}R^{3}}(\mathbf{I}-3\mathbf{u\otimes u})+O\left(  \frac
{1}{R^{2}}\right)  ,\label{per-187}\\
\nabla^{\prime}\times\overline{\mathbf{G}}_{\mathbf{E}}\mathbf{(r}%
,\mathbf{r^{\prime})}  &  =\frac{1}{4\pi R^{2}}\mathbf{u}\times\mathbf{I}%
+O\left(  \frac{1}{R}\right)  , \label{per-189}%
\end{align}
where $\mathbf{u}={(\mathbf{r}}^{\prime}{-\mathbf{r})}/{R},$ which implies
that:
\begin{align}
&  \underset{\delta\rightarrow0}{\text{lim}}\int_{{S}_{\delta}}\ \mathrm{d}%
s^{\prime}\text{ }\mathbf{n\times E}(\mathbf{r}^{\prime})\cdot\nabla
\times\overline{\mathbf{G}}_{E}(\mathbf{r}^{\prime},\mathbf{r})=-\left[
\mathbf{I-L}_{V_{\delta}}\right]  \cdot\mathbf{E}(\mathbf{r}), \label{per-191}%
\\
&  \underset{\delta\rightarrow0}{\text{lim}}\int_{{S}_{\delta}}\ \mathrm{d}%
s^{\prime}\text{ }\mathbf{n}\times\mathbf{H}(\mathbf{r}^{\prime}%
)\cdot\overline{\mathbf{G}}_{E}(\mathbf{r}^{\prime},\mathbf{r})=-\frac
{1}{k^{2}}\mathbf{L}_{V_{\delta}}\cdot\nabla\times\mathbf{H}(\mathbf{r}),
\label{per-193a}%
\end{align}
and the $\mathbf{L}$-dyadics for\ $V_{\delta}$ of various geometric shapes are
given as follows \cite{Yaghjian:1980}:
\begin{equation}
\mathbf{L}_{V_{\delta}}=\frac{1}{3}\mathbf{I}%
\end{equation}
for a sphere.

Substituting (\ref{per-191}) and (\ref{per-193a}) into (\ref{per-186}), after
some manupulation and also using Amp\`{e}re's law, we have the VIE for the
electric field for $\mathbf{r}\in\Omega$:
\begin{equation}
\mathbf{C\cdot E}(\mathbf{r})=\mathbf{E}^{\text{inc}}(\mathbf{r}%
)-\mathrm{i}\omega\mu(\mathbf{r})\text{ p.v.}\int_{\Omega}\mathrm{d}%
\mathbf{r}^{\prime}\text{ }\mathrm{i}\omega\Delta\epsilon(\mathbf{r}^{\prime
})\mathbf{E}(\mathbf{r}^{\prime})\cdot\overline{\mathbf{G}}_{E}(\mathbf{r}%
^{\prime},\mathbf{r}), \label{per-197}%
\end{equation}
where the coefficient matrix is given by
\begin{equation}
\mathbf{C=I+L}_{V_{\delta}}\cdot\Delta\epsilon(\mathbf{r}). \label{per-199}%
\end{equation}

As mentioned in the introduction, in most numerical implementations of
(\ref{per-197}), the CPV integral is computed by selecting a finite $\delta$
as in (\ref{eqn:numericalpv}), therefore the numerical solution thus obtained
does not satisfy the original VIE (\ref{per-197}) and it will have an
intrinsic error reflected in the correction terms indicated in
(\ref{correction}).

\subsection{Reformulation of the VIE and computing CPV with a finite $\delta$}

\label{sec:newvie}

%In the numerical computation of (), a the p.v. is generally approximated by
%selecting a finite and small $\delta,$ namely,%
%\begin{equation}
%\mathbf{C\cdot E}(\mathbf{r})=\mathbf{E}^{\text{inc}}(\mathbf{r}%
%)-\mathrm{i}\omega\mu(\mathbf{r})\text{ }\int_{\Omega\backslash V_{\delta}%
%}\mathrm{d}\mathbf{r}^{\prime}\text{ }\mathrm{i}\omega\Delta\epsilon
%(\mathbf{r}^{\prime})\mathbf{E}(\mathbf{r}^{\prime})\cdot\overline{\mathbf{G}%
%}_{E}(\mathbf{r}^{\prime},\mathbf{r})+O(\delta). \label{vie001}%
%\end{equation}
%Therefore, the numerical solution is in fact an approximation to () by
%ignoring the $O(\delta)$ - which is a "truncation error" of the IE, implying
%that the accuracy of the numerical solution will be always limited by this
%truncation error on the level of IE no matter how accurately the other terms
%of the IE (\ref{vie001}) is evaluated. In this section, we will re-derive the
%VIE (\ref{per-197}) by using vector-scalar potentials for the electric field
%and give the explicit formula for the truncation error in (\ref{vie001}).

In this section, we will derive a VIE of the E-field where the CPV in
(\ref{per-199}) can be computed with a finite exclusion volume together with
some correction terms. Based on the Helmholtz decomposition, the electric
field $\mathbf{E(r)}$ can be expressed as follows \cite{wcai:2013},%
\[
\mathbf{E}=-\mathrm{i}\omega\mathbf{A}-\nabla V,
\]
according to the Lorentz gauge \cite{stratton:1941},
\begin{equation}
\nabla\cdot\mathbf{A}=-\mathrm{i}\omega\epsilon\mu V, \label{eqn1-45}%
\end{equation}

so we have a vector potential representation for the electric field,
%\begin{equation}
%\mathbf{J}(\mathbf{r})=\mathbf{J}_{e}(\mathbf{r})+\mathbf{J}_{\text{eq}%
%}(\mathbf{r}), \label{per-178}%
%\end{equation}%
\begin{equation}
\mathbf{E}=-\mathrm{i}\omega\mathbf{A}+\frac{1}{\mathrm{i}\omega\epsilon\mu
}\nabla(\nabla\cdot\mathbf{A})=-\mathrm{i}\omega\left[  \overline{\mathbf{I}%
}\text{ }\mathbf{+}\text{ }\frac{1}{k^{2}}\nabla\nabla\right]  \mathbf{A}.
\label{eqn1-49}%
\end{equation}

On the other hand, it can be shown that the potential $\mathbf{A}$ satisfies
the Helmholtz equation componentwise \cite{wcai:2013}
\begin{equation}
\nabla^{2}\mathbf{A}+k^{2}\mathbf{A}=-\mu\mathbf{J}. \label{eqn:1-47}%
\end{equation}

Thus the solution $\mathbf{A}$ of Eq. (\ref{eqn:1-47}) can be rewritten as an
integral representation:
\begin{align}
\label{per-180}\mathbf{A}  &  =\mu\int_{\mathbb{R}^{3}}\mathrm{d}%
\mathbf{r}^{\prime}\mathbf{J}(\mathbf{r}^{\prime})g(\mathbf{r},\mathbf{r}%
^{\prime})\nonumber\\
&  =\mu\int_{\mathbb{R}^{3}\backslash\Omega}\mathbf{J}_{e}(\mathbf{r}^{\prime
})g(\mathbf{r},\mathbf{r}^{\prime})\mathrm{d}\mathbf{r}^{\prime}+\mu
\int_{\Omega}\mathbf{J}_{\mathrm{eq}}(\mathbf{r}^{\prime})g(\mathbf{r}%
,\mathbf{r}^{\prime})\mathrm{d}\mathbf{r}^{\prime}\nonumber\\
&  =\mu\int_{\mathbb{R}^{3}\backslash\Omega}\mathbf{J}_{e}(\mathbf{r}^{\prime
})g(\mathbf{r},\mathbf{r}^{\prime})\mathrm{d}\mathbf{r}^{\prime}+\mu
\int_{\Omega}\mathrm{d}\mathbf{r}^{\prime}\mathrm{i}\omega\Delta
\epsilon(\mathbf{r}^{\prime})\mathbf{E}(\mathbf{r}^{\prime})g(\mathbf{r}%
,\mathbf{r}^{\prime}),
\end{align}
where the second equality on the right hand side of Eq. (\ref{per-180}) is due
to the assumption that $\mathrm{supp}(\mathbf{J}_{e}(\mathbf{r}))\cap
\Omega=\emptyset.$

For the first integral in Eq. (\ref{per-180}), it is well-defined if
$\mathbf{r}\in\Omega$ and when we plug it into Eq. (\ref{eqn1-49}), it yields
the incident wave $\mathbf{E}^{\mathrm{inc}}(\mathbf{r})$ according to
relation (\ref{per-173}). For the second integral over $\Omega$ , we split it
as follows
\begin{equation}
\mu\int_{\Omega}\mathrm{d}\mathbf{r}^{\prime}\mathrm{i}\omega\Delta
\epsilon(\mathbf{r}^{\prime})\mathbf{E}(\mathbf{r}^{\prime})g(\mathbf{r}%
,\mathbf{r}^{\prime})=\mu\left(  \int_{\Omega\backslash V_{\delta}}+\int_{V_{\delta}}\right)
\mathrm{i}\omega\Delta\epsilon(\mathbf{r}^{\prime})\mathbf{E}(\mathbf{r}%
^{\prime})g(\mathbf{r},\mathbf{r}^{\prime}),\nonumber
\end{equation}
and along with the first integral term in Eq. (\ref{per-180}), it follows from
Eq. (\ref{eqn1-49}) that%
\begin{align}
\mathbf{E}  &  =\mathbf{E}^{\mathrm{inc}}(\mathbf{r})-\mathrm{i}\omega
\int_{\Omega\backslash V_{\delta}}\mathrm{i}\omega\Delta\epsilon
(\mathbf{r}^{\prime})\mathbf{E}(\mathbf{r}^{\prime})\left[  \overline
{\mathbf{I}}\mathbf{+}\frac{1}{k^{2}}\nabla\nabla\right]  g(\mathbf{r}%
,\mathbf{r}^{\prime})\nonumber\label{per-207}\\
&  -\mathrm{i}\omega\left[  \overline{\mathbf{I}}\mathbf{+}\frac{1}{k^{2}%
}\nabla\nabla\right]  \int_{V_{\delta}}\mathrm{d}\mathbf{r}^{\prime}%
\mathrm{i}\omega\Delta\epsilon(\mathbf{r}^{\prime})\mathbf{E}(\mathbf{r}%
^{\prime})g(\mathbf{r},\mathbf{r}^{\prime})\nonumber\\
&  =\mathbf{E}^{\mathrm{inc}}(\mathbf{r})-\mathrm{i}\omega\int_{\Omega
\backslash V_{\delta}}\mathrm{i}\omega\Delta\epsilon(\mathbf{r}^{\prime
})\overline{\mathbf{G}}_{\mathbf{E}}\mathbf{(r},\mathbf{r^{\prime})\cdot
E}(\mathbf{r}^{\prime})\nonumber\\
&  -\mathrm{i}\omega\left[  \overline{\mathbf{I}}\mathbf{+}\frac{1}{k^{2}%
}\nabla\nabla\right]  \int_{V_{\delta}}\mathrm{d}\mathbf{r}^{\prime}%
\mathrm{i}\omega\Delta\epsilon(\mathbf{r}^{\prime})\mathbf{E}(\mathbf{r}%
^{\prime})g(\mathbf{r},\mathbf{r}^{\prime}).
\end{align}
Next, we isolate the singular part from $g(\mathbf{r},\mathbf{r}^{\prime})$ as
follows
\begin{equation}
g(\mathbf{r},\mathbf{r^{\prime}})=g_{0}(\mathbf{r},\mathbf{r^{\prime}%
})+\widetilde{g}(\mathbf{r},\mathbf{r^{\prime}}), \label{sigsplit}%
\end{equation}
where%
\begin{equation}
g_{0}(\mathbf{r},\mathbf{r^{\prime}})=\frac{1}{4\pi|\mathbf{r}%
-\mathbf{r^{\prime}|}}, \label{per-211}%
\end{equation}
and use the fact that \cite{Lee:1980,Yaghjian:1980}
\begin{equation}
\nabla\nabla\int_{V_{\delta}}\mathrm{d}\mathbf{r}^{\prime}\frac{1}%
{4\pi|\mathbf{r}-\mathbf{r^{\prime}|}}=-\int_{\partial V_{\delta}}%
\mathrm{d}s^{\prime}\frac{\left(  \mathbf{r}-\mathbf{r^{\prime}}\right)
\mathbf{u}_{n}(\mathbf{r^{\prime}})}{4\pi|\mathbf{r}-\mathbf{r^{\prime}|}^{3}%
}=-\mathbf{L}_{V_{\delta}}, \label{per-213}%
\end{equation}
we can compute the following integral in the following manner%
\begin{align}
&  \nabla\nabla\int_{V_{\delta}}\mathrm{d}\mathbf{r}^{\prime}\Delta
\epsilon(\mathbf{r}^{\prime})\mathbf{E}(\mathbf{r}^{\prime})g_{0}%
(\mathbf{r},\mathbf{r}^{\prime})\nonumber\\
&  =\nabla\nabla\int_{V_{\delta}}\mathrm{d}\mathbf{r}^{\prime}\frac{1}%
{4\pi|\mathbf{r}-\mathbf{r^{\prime}|}}\Delta\epsilon(\mathbf{r})\mathbf{E}%
(\mathbf{r})+\int_{V_{\delta}}\mathrm{d}\mathbf{r}^{\prime}\nabla\nabla
g_{0}(\mathbf{r},\mathbf{r}^{\prime})\left[  \Delta\epsilon(\mathbf{r}%
^{\prime})\mathbf{E}(\mathbf{r}^{\prime})-\Delta\epsilon(\mathbf{r}%
)\mathbf{E}(\mathbf{r})\right] \nonumber\\
&  =-\mathbf{L}_{V_{\delta}}\Delta\epsilon(\mathbf{r})\mathbf{E}%
(\mathbf{r})+\int_{V_{\delta}}\mathrm{d}\mathbf{r}^{\prime}\nabla\nabla
g_{0}(\mathbf{r},\mathbf{r}^{\prime})\left[  \Delta\epsilon(\mathbf{r}%
^{\prime})\mathbf{E}(\mathbf{r}^{\prime})-\Delta\epsilon(\mathbf{r}%
)\mathbf{E}(\mathbf{r})\right]  , \label{per-215}%
\end{align}
where the second integral has a removable singularity $O\left(  {\frac
{1}{|\mathbf{r}-\mathbf{r^{\prime}|}^{2}}}\right)  $ through the use of the
spherical coordinate transform centered at $\mathbf{r}$, provided the function
$\Delta\epsilon(\mathbf{r})\mathbf{E}(\mathbf{r})$ is a differentiable inside
$V_{\delta}$, which we assume to be.

Using Eq. (\ref{per-215}) and noting that $\tilde{g}=g-g_{0}$ is a regular
function, Eq. (\ref{per-207}) becomes%
\begin{align}
\mathbf{C\cdot E}  &  =\mathbf{E}^{\text{inc}}(\mathbf{r})-i\omega\mu
\int_{\Omega\backslash V_{\delta}}i\omega\Delta\epsilon(\mathbf{r}^{\prime
})\overline{\mathbf{G}}_{\mathbf{E}}\mathbf{(r},\mathbf{r^{\prime})\cdot
E}(\mathbf{r}^{\prime})\nonumber\\
&  +\omega^{2}\mu\mathrm{\ }\int_{V_{\delta}}\mathrm{d}\mathbf{r}^{\prime
}\Delta\epsilon(\mathbf{r}^{\prime})\mathbf{E}(\mathbf{r}^{\prime
})g(\mathbf{r},\mathbf{r}^{\prime})\nonumber\\
&  +\frac{\omega^{2}}{k^{2}}\mu\mathrm{\ }\int_{V_{\delta}}\mathrm{d}%
\mathbf{r}^{\prime}\Delta\epsilon(\mathbf{r}^{\prime})\nabla\nabla
\widetilde{g}(\mathbf{r},\mathbf{r}^{\prime})\cdot\mathbf{E}(\mathbf{r}%
^{\prime})\nonumber\\
&  +\frac{\omega^{2}}{k^{2}}\mu\int_{V_{\delta}}\mathrm{d}\mathbf{r}^{\prime
}\nabla\nabla g_{0}(\mathbf{r},\mathbf{r}^{\prime})\left[  \Delta
\epsilon(\mathbf{r}^{\prime})\mathbf{E}(\mathbf{r}^{\prime})-\Delta
\epsilon(\mathbf{r})\mathbf{E}(\mathbf{r})\right]  , \label{eqn:per-217}%
\end{align}
with the same coefficient $\mathbf{C}$ being defined in Eq. (\ref{per-199}).

The VIE in (\ref{eqn:per-217}) is similar to those obtained by Fikioris
\cite{Fik:1965}, however, our derivation is based on a splitting of Green's
function in (\ref{sigsplit}) and the identity for the $\mathbf{L}_{V_{\delta}%
}$ in (\ref{per-213}). A comparison study between CPV formulation
(\ref{per-197}) and finite exclusion volume formulation (\ref{eqn:per-217})
can be found in \cite{Wang:1982}. Now expression (\ref{eqn:per-217}) holds for
any finite $\delta>0$ as long as $V_{\delta}\subset\Omega$, and all
integration terms involved on the right hand side are well-defined provided
that $\Delta\epsilon(\mathbf{r})\mathbf{E}(\mathbf{r})$ is H\"{o}lder
continuous. We can see that the last three integration terms can be
understood as the correction terms for computing the Cauchy principal value
with a finite-sized exclusion volume $V_{\delta}$. It should be noted that
these integrals are all weakly singular integrals whose singularities can be
removed by the simple spherical coordinate transform. In particular, we can
easily estimate their magnitude in terms of $\delta.$ Namely,%
\begin{equation}
|\int_{V_{\delta}}\mathrm{d}\mathbf{r}^{\prime}\Delta\epsilon(\mathbf{r}%
^{\prime})\mathbf{E}(\mathbf{r}^{\prime})g(\mathbf{r},\mathbf{r}^{\prime
})|\leq C_{1}||\Delta\epsilon\mathbf{E||}_{\infty}\delta^{2},\label{per-219}%
\end{equation}%
\begin{equation}
|\int_{V_{\delta}}\mathrm{d}\mathbf{r}^{\prime}\Delta\epsilon(\mathbf{r}%
^{\prime})\nabla\nabla\widetilde{g}(\mathbf{r},\mathbf{r}^{\prime}%
)\cdot\mathbf{E}(\mathbf{r}^{\prime})|\leq C_{2}||\Delta\epsilon
\mathbf{E||}_{\infty}\delta^{2},\label{per-221}%
\end{equation}
and
\begin{equation}
|\int_{V_{\delta}}\mathrm{d}\mathbf{r}^{\prime}\nabla\nabla g_{0}%
(\mathbf{r},\mathbf{r}^{\prime})\left[  \Delta\epsilon(\mathbf{r}^{\prime
})\mathbf{E}(\mathbf{r}^{\prime})-\Delta\epsilon(\mathbf{r})\mathbf{E}%
(\mathbf{r})\right]  |\leq C_{3}||\Delta\epsilon\mathbf{E||}_{1,\infty}%
\delta,\label{per-223}%
\end{equation}
where $C_{1}$, $C2$, $C_{3}$ are constants, and $\Vert\cdot\Vert_{\infty}$ and
$\Vert\cdot\Vert_{1,\infty}$ represent the $L^{\infty}$ norms of a function
and its first derivative, respectively.

\begin{rem}
\noindent\ Eqs. (\ref{per-219})-(\ref{per-223}) explicitly indicate the
accuracy of approximating the Cauchy principal value (\ref{cpv}) by the
integral (\ref{eqn:numericalpv}) with a finite $\delta>0$, i.e., the
truncation error is of the order $O(\delta)$. Hence the numerical solution of
the VIE will have this $O(\delta)$ truncation error in general regardless of
the integration quadratures used if terms in Eqs. (\ref{per-219}%
)-(\ref{per-223}) are not included.

However, For the case that $V_{\delta}$ is a ball of radius $\delta$ centered
at $\mathbf{r}$, one can obtain a better estimate in Eq.(\ref{per-223}) due to
the anti-symmetry of the singular term $\nabla\nabla g_{0}(\mathbf{r}%
,\mathbf{r}^{\prime})$ in the spherical coordinates, i.e.
\begin{equation}
|\int_{V_{\delta}}\mathrm{d}\mathbf{r}^{\prime}\nabla\nabla g_{0}%
(\mathbf{r},\mathbf{r}^{\prime})\left[  \Delta\epsilon(\mathbf{r}^{\prime
})\mathbf{E}(\mathbf{r}^{\prime})-\Delta\epsilon(\mathbf{r})\mathbf{E}%
(\mathbf{r})\right]  |\leq C_{4}||\Delta\epsilon\mathbf{E||}_{2,\infty}%
\delta^{2}. \label{eqn:est}%
\end{equation}

\end{rem}

\section{Numerical Methods}

\label{sec:methods}

\subsection{Nystr\"{o}m collocation method}

Now we will use the Nystr\"{o}m collocation method to solve Eq.
(\ref{eqn:per-217}). First, the computational domain (scatter) $\Omega$ is
divided into $N$ cubic elements $\Omega_{i}$ with length $a_{i},i=1,2,...,N$.
On each element $\Omega_{i}$, we assign $M$ Gauss nodes on which $M$ scalar
Lagrange basis functions $\phi_{ij},j=1,2,3,....M$ are defined and vanish
outside $\Omega_{i}$. Then, we can write the solution as
\begin{equation}
\mathbf{E}({}\mathbf{r})=\sum_{i=1}^{N}\sum_{j=1}^{M}\mathbf{c}_{ij}\phi
_{ij}(\mathbf{r}),\text{ \ \ \ \ }\mathbf{r\in}\Omega_{i},
\label{eqn:discretization}%
\end{equation}
with $\mathbf{c}_{ij}$ are the unknown vectorial coefficients. Inserting Eq.
(\ref{eqn:discretization}) into Eq. (\ref{eqn:per-217}),
%\begin{equation}\nonumber
%	\phi_{ij}({\bf r}_{\tilde{i}\tilde{j}})=1 \quad\text{ only when } i=\tilde{i} \text{ and } j=\tilde{j}
%\end{equation}
we obtain the following equations for $\mathbf{c}_{ij}$:
%{\small\begin{eqnarray}\nonumber
%	\mathbf{C} \cdot \sum_{i=1}^N\sum_{j=1}^M\mathbf{c}_{ij}\phi_{ij}({\bf r}_{ij})&=&\mathbf{E}^{\rm inc}({\bf r}_{ij})+\omega^2\mu\sum_{i=1}^N\sum_{j=1}^M\left[\int_{\Omega\backslash V_{\delta_{ij}}}\mathrm{d}{\bf r}'\Delta\epsilon({\bf r}')\bar{\bf G}_{\bf E}({\bf r}_{ij}, {\bf r}')\phi_{ij}({\bf r}')\right]{\bf c}_{ij}\\\nonumber
%	&+&\omega^2\mu\sum_{i=1}^N\sum_{j=1}^M\left[\int_{V_{\delta_{ij}}}\mathrm{d}{\bf r}'\Delta\epsilon(\mathbf{r}^{\prime})g(\mathbf{r}_{ij},\mathbf{r}^{\prime})\phi_{ij}({\bf r}')\right]{\bf c}_{ij}\\\nonumber
%	&+&\frac{\omega^2\mu}{k^2_L}\sum_{i=1}^N\sum_{j=1}^M\left[\int_{V_{\delta_{ij}}}\mathrm{d}{\bf r}'\Delta\epsilon(\mathbf{r}^{\prime})\nabla\nabla
%\widetilde{g}(\mathbf{r}_{ij},\mathbf{r}^{\prime})\phi_{ij}({\bf r}')\right]\cdot{\bf c}_{ij}\\
%	&+&\frac{\omega^{2}\mu}{k_{L}^{2}}{\int_{V_{\delta_{ij}}}\mathrm{d}\mathbf{r}%
%^{\prime}\nabla\nabla g_{0}(\mathbf{r}_{ij},\mathbf{r}^{\prime})\left[
%\Delta\epsilon(\mathbf{r}^{\prime})\phi_{ij}(\mathbf{r}^{\prime}%
%)-\Delta\epsilon(\mathbf{r}_{ij})\phi_{ij}(\mathbf{r}_{ij})\right]}\cdot{\bf c}_{ij}
%\end{eqnarray}}
{\small
\begin{align}
\mathbf{C}\cdot\mathbf{c}_{ij}  &  =\mathbf{E}_{ij}^{\mathrm{inc}}+\omega
^{2}\mu\sum_{n=1}^{N}\sum_{m=1}^{M}\left[  \int_{\Omega_{n}\backslash
V_{\delta_{ij}}}\mathrm{d}\mathbf{r}^{\prime}\Delta\epsilon(\mathbf{r}%
^{\prime})\bar{\mathbf{G}}_{\mathbf{E}}(\mathbf{r}_{ij},\mathbf{r}^{\prime
})\phi_{nm}(\mathbf{r}^{\prime})\right]  \mathbf{c}_{nm}%
\nonumber\label{eqn:discre}\\
&  +\omega^{2}\mu\sum_{m=1}^{M}\left[  \int_{V_{\delta_{ij}}}\mathrm{d}%
\mathbf{r}^{\prime}\Delta\epsilon(\mathbf{r}^{\prime})g(\mathbf{r}%
_{ij},\mathbf{r}^{\prime})\phi_{im}(\mathbf{r}^{\prime})\right]
\mathbf{c}_{im}\nonumber\\
&  +\frac{\omega^{2}\mu}{k^{2}}\sum_{m=1}^{M}\left[  \int_{V_{\delta_{ij}%
}}\mathrm{d}\mathbf{r}^{\prime}\Delta\epsilon(\mathbf{r}^{\prime})\nabla
\nabla\widetilde{g}(\mathbf{r}_{ij},\mathbf{r}^{\prime})\phi_{im}%
(\mathbf{r}^{\prime})\right]  \cdot\mathbf{c}_{im}\nonumber\\
&  +\frac{\omega^{2}\mu}{k^{2}}\sum_{m=1}^{M}{\int_{V_{\delta_{ij}}%
}\mathrm{d}\mathbf{r}^{\prime}\nabla^{2}g_{0}(\mathbf{r}_{ij},\mathbf{r}%
^{\prime})\left[  \Delta\epsilon(\mathbf{r}^{\prime})\phi_{im}(\mathbf{r}%
^{\prime})-\Delta\epsilon_{ij}\phi_{im}(\mathbf{r}_{ij})\right]  }%
\cdot\mathbf{c}_{im}.
\end{align}
}

When $n\neq i$, we have $V_{\delta_{ij}}\notin\Omega$, then in Eq.
(\ref{eqn:discre}) the integral
\begin{equation}
\int_{\Omega_{n}\backslash V_{\delta_{ij}}}\mathrm{d}\mathbf{r}^{\prime}%
\Delta\epsilon(\mathbf{r}^{\prime})\bar{\mathbf{G}}_{\mathbf{E}}%
(\mathbf{r}_{ij},\mathbf{r}^{\prime})\phi_{nm}(\mathbf{r}^{\prime})
\end{equation}
is regular on the whole cube $\Omega_{n}$ and hence it can be evaluated by the
regular Gauss quadratures, i.e.,
\begin{equation}
\int_{\Omega_{n}}\mathrm{d}\mathbf{r}^{\prime}\Delta\epsilon(\mathbf{r}%
^{\prime})\bar{\mathbf{G}}_{\mathbf{E}}(\mathbf{r}_{ij},\mathbf{r}^{\prime
})\phi_{nm}(\mathbf{r}^{\prime})=\left(  \frac{a_{n}}{2}\right)  ^{3}%
\sum_{m=1}^{M}\Delta\epsilon_{nm}\bar{\mathbf{G}}_{\mathbf{E}}(\mathbf{r}%
_{ij},\mathbf{r}_{nm})\omega_{m}^{s},
\end{equation}
with $\omega_{m}^{s}$ being the standard Gauss weights in 3-D, which are
obtained from the tensor product of the Gauss weights in the 1-D domain
$[-1,1]$.
%\begin{eqnarray}\nonumber
%	\sum_{n=1}^N\sum_{m=1}^M\left[\int_{\Omega_n\backslash V_{\delta_{\tilde{i}\tilde{j}}}}\mathrm{d}{\bf r}'\Delta\epsilon\bar{\bf G}_{\bf E}({\bf r}_{\tilde{i}\tilde{j}}, {\bf r}')\phi_{ij}({\bf r}')\right]{\bf c}_{ij}&=&\sum_{j=1}^M\left[\int_{\Omega_{\tilde{i}}\backslash V_{\delta_{\tilde{i}\tilde{j}}}}\mathrm{d}{\bf r}'\Delta\epsilon\bar{\bf G}_{\bf E}({\bf r}_{\tilde{i}\tilde{j}}, {\bf r}')\phi_{\tilde{i}j}({\bf r}')\right]{\bf c}_{\tilde{i}j}\\\nonumber
%&+& \sum_{i=1, i\neq\tilde{i}}^N\sum_{j=1}^M\left[\int_{\Omega_i}\mathrm{d}{\bf r}'\Delta\epsilon\bar{\bf G}_{\bf E}({\bf r}_{\tilde{i}\tilde{j}}, {\bf r}')\phi_{ij}({\bf r}')\right]{\bf c}_{ij}
%\end{eqnarray}
%where $\Omega_i$ is the $i$-th element.

When $n=i$, although the singularity $\mathbf{r}_{ij}$ is excluded from the
domain $\Omega_{i}$, the calculation of the integral
\begin{equation}
\int_{\Omega_{i}\backslash V_{\delta_{ij}}}\mathrm{d}\mathbf{r}^{\prime}%
\Delta\epsilon(\mathbf{r}^{\prime})\bar{\mathbf{G}}_{\mathbf{E}}%
(\mathbf{r}_{ij},\mathbf{r}^{\prime})\phi_{im}(\mathbf{r}^{\prime})
\label{eqn:singular}%
\end{equation}
is still challenging. Next, we present an efficient quadrature formula to
evaluate this integral.

\subsection{Interpolated weights on Gauss nodes for integrals on
$\Omega\backslash V_{\delta}$}

For the sake of generality, we consider the following integral:
\begin{equation}
I_{s}=\int_{\Omega\backslash V_{\delta}}\frac{f(\mathbf{r};\mathbf{r}^{\prime
})h(\mathbf{r};\mathbf{r}^{\prime})}{R^{k}}\mathrm{d}\mathbf{r}^{\prime}%
,\quad\mathbf{r}\in V_{\delta}, \label{eqn:singular-general}%
\end{equation}
where $k=1,2,3$ corresponds to weak, strong, and hyper-singularity of the
integral, respectively. The function $f(\mathbf{r};\mathbf{r}^{\prime})$ is
assumed to be a general smooth and well-defined function, while $h(\mathbf{r}%
;\mathbf{r}^{\prime})$ is a fixed arbitrary function which results from the
directional derivative in the definition of the dyadic Green's functions.

As the function $f(\mathbf{r};\mathbf{r}^{\prime})$ is smooth over the whole
domain $\Omega$, then it can be well approximated by the following simple
interpolation:
\begin{equation}
f(\mathbf{r};\mathbf{r}^{\prime})\approx\sum_{j=1}^{J}f(\mathbf{r}%
;\mathbf{r}_{j})\phi_{j}(\mathbf{r}^{\prime}),\text{ \ }\mathbf{r}_{j}%
\in\Omega, \label{eqn:interpolation}%
\end{equation}
where $\{\mathbf{r}_{j}\}_{j=1}^{J}$ are $J$-nodes in $\Omega$ formed by the
tensor product of Gauss node in $[-1,1]$. Combining Eqs.
(\ref{eqn:singular-general}) and (\ref{eqn:interpolation}) yields
\begin{equation}
\int_{\Omega\backslash V_{\delta}}\frac{f(\mathbf{r};\mathbf{r}^{\prime
})h(\mathbf{r};\mathbf{r}^{\prime})}{R^{k}}d\mathbf{r}^{\prime}\approx
\sum_{j=1}^{J}f(\mathbf{r};\mathbf{r}_{j})\omega_{j}^{e},
\end{equation}
We call $\omega_{j}^{e}$ the effective interpolated weights, or just
\emph{interpolated weights}, which are defined through the integral
\begin{equation}
\omega_{j}^{e}=\int_{\Omega\backslash V_{\delta}}\frac{\phi_{j}(\mathbf{r}%
^{\prime})h(\mathbf{r};\mathbf{r}^{\prime})}{|\mathbf{r}-\mathbf{r}^{\prime
}|^{k}}d\mathbf{r}^{\prime}. \label{eqn:weights}%
\end{equation}
Note that the interpolated weights $\omega_{j}^{e}$ depend on the location of
the singularity $\mathbf{r}$ and relies on accurate calculations of Eq.
(\ref{eqn:weights}), which will be accomplished in two steps: first the domain
$\Omega$ is subdivided into nine cubes with one containing the singularity
$\mathbf{r}$ in its center. For the cube including the singularity, a
straight-forward, brute-force approach involving a large number $N$ ($N\gg J$)
of Gauss points is adopted in local spherical polar coordinates, will be used
to obtain satisfactory accuracy. While for the other cubes, regular tensor
product Gauss quadrature is applied. Details of the computations and resulting
weight tables could be found in \cite{Zinser:2015}.

Fortunately, the computation of weights $\omega_{j}^{e}$ only needs to be
performed once and tabulated for the reference domain, then they can be used
for varies cubic elements. Due to the smoothness of function $f(\mathbf{r}%
,\mathbf{r}^{\prime})$, the number $J$ is relatively small, especially if the
element size is small as in practice, so computation of Eq.
(\ref{eqn:singular}) is now efficient once $\omega_{j}^{e}$ are obtained.

\begin{rem}
As the Cauchy principal value in the VIE, thus the VIE itself, depends on the
specific shape of the exclusion volume, the interpolated weights can only be
used for the shape for which it was calculated, namely, a cubic shape here.
So, for a general element obtained by a affine mapping as in a finite element
triangulation, we will need to isolate the singularity by a cube with the
singularity at its center, then the pre-calculated interpolated weights
defined in (\ref{eqn:weights}) can be used while integral over the rest of the
region outside the cube within the element of a more general shape can be done
with regular Gauss quadratures.
\end{rem}

\subsection{Computation of VIE matrix entries}

In this section, we will show how to compute the matrix entries accurately for
the VIE in the following steps.

\begin{itemize}
\item Step I: calculate interpolated weights on the reference cubic domain.

\bigskip

For Eq. (\ref{eqn:weights}), we take $\Omega=[-1,1]^{3}$ and $V_{\delta
}=B(\mathbf{r}_{j},{\delta})$, where $\delta>0$ is a prescribed small number
and $\mathbf{r}_{j},j=1,2,...M$ coincident with the coordinates of $M$ Gauss
points in the reference domain, and we take $J=M$ in Eq.
(\ref{eqn:interpolation}).

Further, the dyadic Green's function has the form
\begin{align}
\overline{\mathbf{G}}_{\mathbf{E}}  &  =g\mathbf{I}+\frac{\nabla^{2}g}%
{k^{2}}=\frac{e^{-ikR}}{4\pi R}(\mathbf{I}-\mathbf{u}\otimes
\mathbf{u})\nonumber\label{eqn:G-full}\\
&  -\frac{ie^{-ikR}}{4\pi R^{2}k}(\mathbf{I}-3\mathbf{u}%
\otimes\mathbf{u})-\frac{e^{-ikR}}{4\pi R^{3}k^{2}}(\mathbf{I}%
-3\mathbf{u}\otimes\mathbf{u}).
\end{align}
%where $R=|{\bf r}-{\bf r}'|$ and ${\bf u}$ is the unit vector of ${\bf r}-{\bf r}'$, i.e., ${\bf u}=\displaystyle{\frac{{\bf r}-{\bf r}'}{R}}$,
Since the value of function $\mathbf{u}\otimes\mathbf{u}$ is multiple-defined
at $R=0$, in (\ref{eqn:singular-general}) we need to take
\begin{equation}
h(\mathbf{r};\mathbf{r}^{\prime})=\mathbf{u}\otimes\mathbf{u,}%
\end{equation}
and hence Eq. (\ref{eqn:weights}) yields a set of $9$ interpolated weight
matrices. However, due to the symmetry of the matrix $\mathbf{u}%
\otimes\mathbf{u}$, only 6 components need to be calculated. For the identity
matrix $\mathbf{I}$ term in Eq. (\ref{eqn:G-full}), we will also need a set of
scalar interpolation weights by assuming $h(\mathbf{r},\mathbf{r}^{\prime})=1$
in (\ref{eqn:weights}).

Additionally, for the scalar and matrix weights, we need to calculate for
$k=1,2$, and $3$ for weak, strong, and hyper singular integrals, respectively.

To summarize, for each collocation point (also the singularity location)
$\mathbf{r}_{j},j=1,2,...,M$ in an element, scalar weights $\omega_{j,m}$,
$\bar{\omega}_{j,m}$ and $\tilde{\omega}_{j,m},m=1,2,...,M$ are calculated
for weak, strong, hyper-singular integrals, respectively. And the
corresponding weights are denoted as $\Lambda_{j,m}$, $\bar{\Lambda}_{j,m}%
$, and $\tilde{\Lambda}_{j,m}$. These weights only need to be calculated
once and then stored for future use.

\bigskip

\item Step II: We discretize the computational domain into a series of cubic
elements with length $a_{i},i=1,2,...,N$ and assign $M$ Gauss points in each
elements. For the $j$-th Gauss point $\mathbf{r}_{ij}$ in the $i$-th element,
we construct the equation:
\begin{equation}
-\omega^{2}\mu\sum_{n=1}^{N}\sum_{m=1}^{M}\mathbf{A}_{nm}\cdot\mathbf{c}%
_{nm}-\sum_{m=1}^{M}\mathbf{B}_{im}\cdot\mathbf{c}_{im}+\left(  1+\frac{1}%
{3}\Delta\epsilon_{ij}\right)  \mathbf{I}_{3\times3}\cdot\mathbf{c}%
_{ij}=\mathbf{E}_{ij}^{\mathrm{inc}}, \label{eqn:equation}%
\end{equation}
%		The coefficients
%		\begin{equation}
%			{\bf C}=\left(1+\frac{1}{3}\Delta\epsilon_{ij}\right){\bf I}
%		\end{equation}
The matrix $\mathbf{B}$ originates from the correction terms of the Cauchy
principal value
\begin{align}
\mathbf{B}_{im}  &  =\omega^{2}\mu\int_{B(\mathbf{r}_{ij},a_{i}\delta
)}\mathrm{d}\mathbf{r}^{\prime}\Delta\epsilon(\mathbf{r}^{\prime}%
)g(\mathbf{r}_{ij},\mathbf{r}^{\prime})\phi_{im}(\mathbf{r}^{\prime
})\nonumber\\
&  +\frac{\omega^{2}\mu}{k^{2}}\int_{B(\mathbf{r}_{ij},a_{i}\delta)}%
\mathrm{d}\mathbf{r}^{\prime}\Delta\epsilon(\mathbf{r}^{\prime})\nabla
^{2}\tilde{g}(\mathbf{r}_{ij},\mathbf{r}^{\prime})\phi_{im}(\mathbf{r}%
^{\prime})\nonumber\\
&  +\frac{\omega^{2}\mu}{k^{2}}\int_{B(\mathbf{r}_{ij},a_{i}\delta)}%
\mathrm{d}\mathbf{r}^{\prime}\nabla^{2}g_{0}(\mathbf{r}_{ij},\mathbf{r}%
^{\prime})\left[  \Delta\epsilon(\mathbf{r}^{\prime})\phi_{im}(\mathbf{r}%
^{\prime})-\Delta\epsilon_{ij}\phi_{im}(\mathbf{r}_{ij}))\right]  ,
\end{align}
and it can be calculated by standard Gauss quadrature through spherical
coordinates since the Jacobian will eliminate completely the singularity of
the integrands.

When $n=i$, we calculate the integral Eq.(\ref{eqn:singular}) as
\begin{align}
\mathbf{A}_{im}  &  =\frac{1}{4\pi}\left(  \frac{a_{i}}{2}\right)  ^{3}%
\sum_{j=1}^{M}\Delta\epsilon_{im}\left[  \left(  e^{-ikR_{m}}\omega
_{j,m}^{i}-ie^{-ikR_{m}}\bar{\omega}_{j,m}^{i}-e^{-ikR_{m}}\tilde{\omega
}_{j,m}^{i}\right)  \mathbf{I}_{3\times3}\right. \nonumber\\
&  \left.  e^{-ikR_{m}}\Lambda_{j,m}^{i}-ie^{-ikR_{m}}\bar{\Lambda}%
_{j,m}^{i}-e^{-ikR_{m}}\tilde{\Lambda}_{j,m}^{i}\right]  ,
\end{align}
where $R_{m}=|\mathbf{r}_{ij}-\mathbf{r}_{im}|$ and recall the definition of
the interpolated weights Eq. (\ref{eqn:weights}), we have
\begin{equation}%
\begin{array}[c]{lll}%
\omega_{j,m}^{i}=\left(  \frac{2}{a_{i}}\right)  \omega_{j,m}, & \bar{\omega
}_{j,m}^{i}=\left(  \frac{2}{a_{i}}\right)  ^{2}\bar{\omega}_{j,m}, &
\tilde{\omega}_{j,m}^{i}=\left(  \frac{2}{a_{i}}\right)  ^{3}\tilde{\omega}%
_{j,m}\\
\label{eqn:interpolated-weights}\Lambda_{m}^{i}=\left(  \frac{2}{a_{i}%
}\right)  \Lambda_{j,m}, & \bar{\Lambda}_{j,m}^{i}=\left(  \frac{2}{a_{i}%
}\right)  ^{2}\bar{\Lambda}_{j,m}, & \tilde{\Lambda}_{j,m}^{i}=\left(
\frac{2}{a_{i}}\right)  ^{3}\tilde{\Lambda}_{j,m}%
\end{array}
,
\end{equation}
when $n\neq i$, we have
\begin{equation}
\mathbf{A}_{nm}=\left(  \frac{a_{n}}{2}\right)  ^{3}\sum_{j=1}^{M}%
\Delta\epsilon_{nm}\bar{\mathbf{G}}_{\mathbf{E}}(\mathbf{r}_{ij}%
,\mathbf{r}_{nm})\omega_{j}^{s}.
\end{equation}

%Map the reference domain to the physical element with length $a$ excluding a ball with radius $a\delta$.

\item Step III: Eq. (\ref{eqn:equation}) for all the $N\times M$ Gauss points
can be assembled as the following linear algebraic equation system
\begin{equation}
\mathbf{V}\cdot\vec{\mathbf{c}}=\left[
\begin{array}
[c]{ccc}%
\mathbf{V}_{xx} & \mathbf{V}_{xy} & \mathbf{V}_{xz}\\
\mathbf{V}_{yx} & \mathbf{V}_{yy} & \mathbf{V}_{yz}\\
\mathbf{V}_{zx} & \mathbf{V}_{zy} & \mathbf{V}_{zz}%
\end{array}
\right]  \cdot\left[
\begin{array}
[c]{c}%
\mathbf{c}_{x}\\
\mathbf{c}_{y}\\
\mathbf{c}_{z}%
\end{array}
\right]  =\left[
\begin{array}
[c]{c}%
\mathbf{E}_{x}^{\mathrm{inc}}\\
\mathbf{E}_{y}^{\mathrm{inc}}\\
\mathbf{E}_{z}^{\mathrm{inc}}%
\end{array}
\right]  . \label{eqn:sys}%
\end{equation}
Based on the properties of the Green's function, the $3NM\times3NM$ matrix
$\mathbf{V}$ is partitioned into nine blocks, each of which is a $NM\times NM$
sub-matrix. The solution of the VIE contains three $NM\times1$ vectors, which
represents the field in $x$, $y$, and $z$ directions. This system is solved by
iteration methods such as the GMRES method.
\end{itemize}

\section{Numerical Results}

\label{sec:results}

In this section we test the accuracy of the interpolated weights on Gauss
nodes, the $\delta$-independence of the solution of the VIE, and convergence
of the $p$-refinement of the Nystr\"{o}m collocation method.

\subsection{Accuracy of the interpolated weights on Gauss nodes}

In Eq. (\ref{eqn:equation}), the calculation of matrix $\mathbf{B}$ from the
correction terms are straightforward, so we will focus on validating the
interpolated weights in computing matrix $\mathbf{A}$. For convenience, we
consider the integral of a real-valued, tensor function
\begin{equation}
\frac{\cos{R}}{R}(\mathbf{I}-\mathbf{u}\otimes\mathbf{u})+\frac{\cos{R}}%
{R^{2}}(\mathbf{I}-3\mathbf{u}\otimes\mathbf{u})+\frac{\cos{R}}{R^{3}%
}(\mathbf{I}-3\mathbf{u}\otimes\mathbf{u}), \label{eqn:sampleintegral}%
\end{equation}
which is similar to the Green's function Eq. (\ref{eqn:G-full}) on the domain
$\Omega\backslash V_{\delta}$. Without loss of generality, we take
$\Omega=[-1,1]^{3}$ and $\mathbf{r}_j$ as the 27 points constructed from the tensor
product of the Gauss points of order 3 in $[-1,1]$. Thus, we have
$j=m=1,2,...,27$ as in Eq.(\ref{eqn:interpolated-weights}) and the integral
yields
\begin{equation}
\mathbf{G}_{j}\approx\sum_{m=1}^{27}\cos(|\mathbf{r}_{m}-\mathbf{r}%
_{j}|)\left[  \left(  \omega_{j,m}+\bar{\omega}_{j,m}+\tilde{\omega}%
_{j,m}\right) \bar{ \mathbf{I}}-\Lambda_{j,m}-3\bar{\Lambda}_{j,m}%
-3\tilde{\Lambda}_{j,m}\right]  , \label{eqn:appxofG}%
\end{equation}
which is a $3\times3$ matrix depending on $\mathbf{r}_{j}$.
%		we test the accuracy of the interpolated Gauss weights by the sample integral
%		\begin{equation}\label{eqn:sampleintegral}
%			{\bf G}=\int_{\Omega\backslash V_{\delta}}\left[\frac{\cos{R}}{R}({\bf I}-{\bf u}\otimes{\bf u})+\frac{\cos{R}}{R^2}({\bf I}-3{\bf u}\otimes{\bf u})+\frac{\cos{R}}{R^3}({\bf I}-3{\bf u}\otimes{\bf u})\right]\mathrm{d}{\bf r}',
%		\end{equation}
%		
%		  and Eq. (\ref{eqn:sampleintegral}) is approximated as

For each $\mathbf{G}_{j}$, we use the direct brute force method introduced in
\cite{Zinser:2015} to obtain the reference solution with the very small
$\delta=10^{-3}$. Then we calculate the integral using the interpolated
weights as in Eq.(\ref{eqn:appxofG}) with different values of $\delta$.
According to the previous analysis, the differences with the reference
solution should decay in the order $O(\delta^{2})$ as $\delta$ decreases,
which will be checked in the following.

We classify the 27 sets of weights into four categories, based on the position
of the singularity $\mathbf{r}_j$, as the ones near the corner, edge, face,
and center of the cube. The matrix $\mathbf{G}_j$ is symmetric, so we only
check the three diagonal entries ($g_{11}$, $g_{22}$, and $g_{33}$) and the
three  upper diagonal entries (${g_{12}}$, $g_{13}$, and $g_{23}$).

Table \ref{table:center} presents the numerical results when the singularity
$\mathbf{r}_{j}$ is located in the center of the cube, in which case
$g_{11}=g_{22}=g_{33}$ and the off-diagonal entries are all zeros. The values
of $\delta$ are taken as $0.1$, $0.05$,$0.025$, and $0.0125$ while the
reference solution is $g_{11}=4.027477$. 

\begin{table}[ptb]
\caption{Convergence of the integral as the singularity is at center.
$g_{11}=g_{22}=g_{33}$ and $g_{12}=g_{13}=g_{23}=0$.}%
\label{table:center}%
\centering
\begin{tabular}
[c]{l|l|l|l|l}\hline\hline
& $\delta= 0.1$ & $\delta= 0.05$ & $\delta= 0.025$ & $\delta= 0.0125$\\\hline
$g_{11}$ & 3.985701 & 4.017024 & 4.024872 & 4.026835\\
error & -4.1784E-2 & -1.0461E-2 & -2.613E-3 & -6.5E-4\\
order & - & 2 & 2 & 2\\\hline\hline
\end{tabular}
\end{table}

\begin{table}[ptb]
\caption{Convergence of the integral as singularity is at a corner.
Reference solution $g_{11}=g_{22}=g_{33}=0.982526$ and $g_{12}=g_{13}%
=g_{23}=-0.998097$.}%
\label{table:corner}%
\centering
\begin{tabular}
[c]{l|l|l|l|l}\hline\hline
& $\delta= 0.1$ & $\delta= 0.05$ & $\delta= 0.025$ & $\delta= 0.0125$\\\hline
$g_{11}$ & 0.940714 & 0.972063 & 0.979913 & 0.981876\\
error & -2.80425E-1 & -7.3424E-2 & -1.8102E-2 & -3.983E-3\\
order & - & 1.93 & 2 & 2.1\\\hline\hline
$g_{12}$ & -0.998084 & -0.998094 & -0.998097 & -0.998097\\
error & 1.3E-5 & 3.0E-6 & 0 & 0\\
order & - & 2.1 & - & -\\\hline\hline
\end{tabular}
\end{table}

\begin{table}[ptb]
\caption{Convergence of the integral as singularity is at an edge.
Reference solutions $g_{11}=-1.515302$, $g_{22}=g_{33}=3.39234$,
$g_{23}=-1.579086$ and $g_{12}=g_{13}=0$.}%
\label{table:edge}%
\centering
\begin{tabular}
[c]{l|l|l|l|l}\hline\hline
& $\delta= 0.1$ & $\delta= 0.05$ & $\delta= 0.025$ & $\delta= 0.0125$\\\hline
$g_{11}$ & -1.559532 & -1.526351 & -1.518059 & -1.515987\\
error & -4.423E-2 & 1.1049E-2 & 2.757E-3 & -6.85E-4\\
order & - & 2 & 2 & 2\\\hline\hline
$g_{22}$ & 3.35175 & 3.38217 & 3.389798 & 3.391707\\
error & -4.059E-2 & -1.1017E-2 & 2.542E-3 & -6.33E-4\\
order & - & 1.9 & 2.1 & 2\\\hline\hline
$g_{23}$ & -1.579072 & -1.579082 & -1.579085 & -1.579085\\
error & 1.4E-5 & 4.0E-6 & 1.0E-6 & 1.0E-6\\
order & - & 1.8 & 2 & 0\\\hline\hline
\end{tabular}
\end{table}

\begin{table}[ptb]
\caption{Convergence of the integral as singularity is at a face.
Reference solutions $g_{11}=0.877428$, $g_{22}=0$, $g_{33}=6.494784$, and
$g_{12}=g_{13}=g_{23}=0$.}%
\label{table:face}%
\centering
\begin{tabular}
[c]{l|l|l|l|l}\hline\hline
& $\delta= 0.1$ & $\delta= 0.05$ & $\delta= 0.025$ & $\delta= 0.0125$\\\hline
$g_{11}$ & 0.83442 & 0.866672 & 0.874742 & 0.87676\\
error & -4.3008E-2 & -1.10756E-2 & -2.686E-3 & -6.68E-4\\
order & - & 1.9 & 2 & 2\\\hline\hline
$g_{33}$ & 6.455419 & 6.484909 & 6.492315 & 6.49417\\
error & -3.9365E-2 & -9.875E-3 & -2.469E-3 & -6.14E-4\\
order & - & 2 & 2 & 2\\\hline\hline
\end{tabular}
\end{table}

In a similar fashion, Tables \ref{table:corner}-\ref{table:face} show the
accuracies when the singularity is located near the corner, edge, and face of
the cube, respectively. Comparing to the reference solutions, the expected
second-order decay with respect to the $\delta$ is confirmed.

\bigskip

\subsection{$\delta$-independence of the VIE solution}

\begin{figure}[ptb]
\begin{center}%
\begin{tabular}
[c]{cc}%
\includegraphics[width=0.5\textwidth]{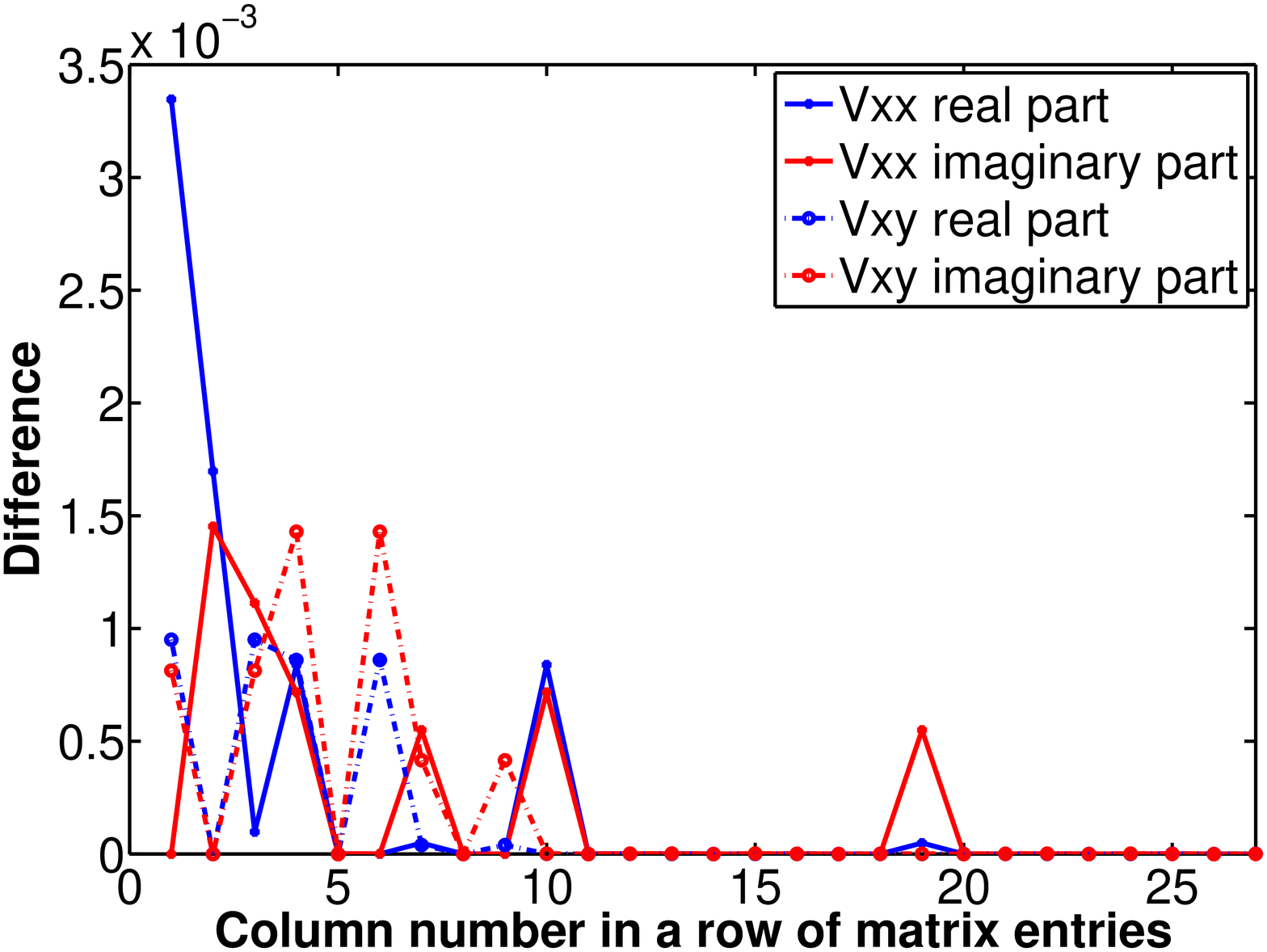} &
\includegraphics[width=0.5\textwidth]{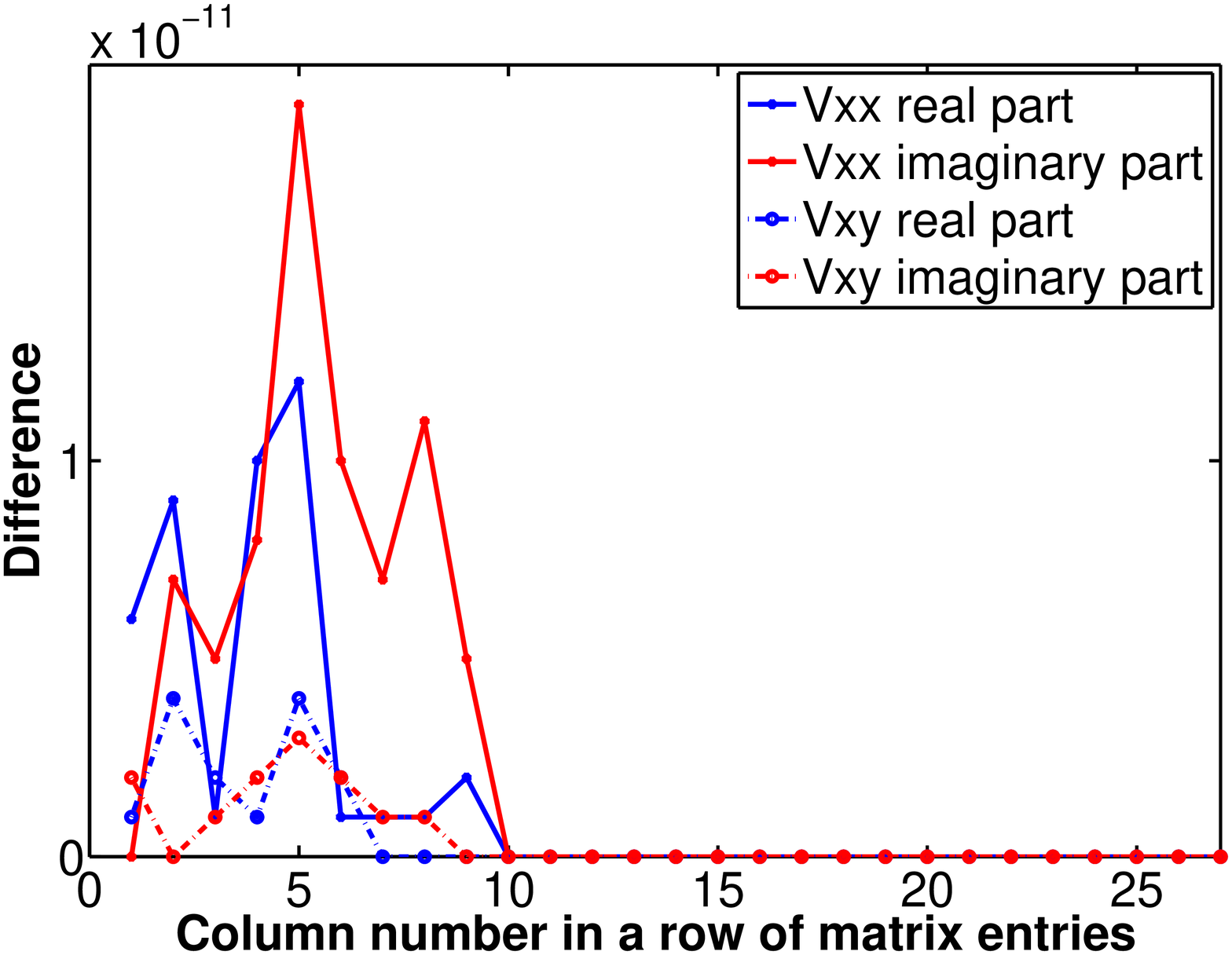}\\
(a) & (b)
\end{tabular}
\end{center}
\caption{Differences of matrix entries with $\delta$=0.1 and $\delta=0.001$.
(a): without correction terms; (b) with correction terms.}%
\label{fig:entryerror}%
\end{figure}

Equation (\ref{eqn:per-217}) provides a formulation from which the solution of
the VIE will be independent of the choice of a specific $\delta$. In the
following tests, we take $\mu=1$, $\Delta\epsilon=4$, $\omega=1$ and solve the VIE. The
computation domain is taken as $[-\pi/2,\pi/2]^{3}$, while the incident wave
is
\begin{equation}
\mathbf{E}_{x}^{\mathrm{inc}}=e^{ik(-y+0.5z)},\quad\mathbf{E}_{y}%
^{\mathrm{inc}}=\mathbf{E}_{z}^{\mathrm{inc}}=0.\nonumber
\end{equation}
We first check the $\delta$-independence of the matrix entries in Eq.
(\ref{eqn:sys}). Figure \ref{fig:entryerror} displays the differences of one
row of entries in the matrix $\mathbf{V}$ between the choices of $\delta=0.1$
and $\delta=0.001$, in which the solid lines are for the entries from a
diagonal block ($V_{xx}$) and dashed lines are for the entries from an
off-diagonal block ($V_{xy}$). The blue curves are for real while red curves
are for imaginary parts. From Fig. \ref{fig:entryerror}(a) we can see that
the differences between entries in the corresponding positions can be as large
as $3.0\times10^{-3}$ when the correction terms are not included. In contrast,
the corresponding differences are reduced to below $2\times10^{-11}$ when the
corrections are, hence the matrix entries are $\delta$-independent.

Next we check the $\delta$-dependence of the overall solution of the VIE. The
solution of VIE with a very small $\delta=0.001$ is taken as the reference
solution. Then various choices of $\delta$ are taken and the differences of
the corresponding solutions with the reference solution are calculated. The
differences are measured in the $L^{\infty}$ norm for the three components
$E_{x}$, $E_{y}$, $E_{z}$ and they are listed in Tables \ref{table:sol1}%
-\ref{table:sol2}. \begin{table}[th]
\caption{Comparison of solutions of the VIE without the correction terms}%
\label{table:sol1}%
\centering
\begin{tabular}
[c]{l|l|l|l|l}\hline\hline
& $\delta= 0.1$ & $\delta= 0.05$ & $\delta= 0.025$ & $\delta= 0.0125$\\\hline
$\|E_{x}-E_{x}^{\mathrm{ref}}\|_{L^{\infty}}$ & 3.360E-3 & 8.264E-4 &
2.044E-4 & 5.039E-5\\
$\|E_{y}-E_{y}^{\mathrm{ref}}\|_{L^{\infty}}$ & 1.476E-3 & 3.696E-4 &
9.258E-5 & 2.358E-5\\
$\|E_{z}-E_{z}^{\mathrm{ref}}\|_{L^{\infty}}$ & 2.533E-3 & 6.387E-4 &
1.597E-4 & 3.969E-5\\\hline\hline
\end{tabular}
\end{table}\begin{table}[thpth]
\caption{Comparison of solutions of the VIE with the correction terms}%
\label{table:sol2}%
\centering
\begin{tabular}
[c]{l|l|l|l|l}\hline\hline
& $\delta= 0.1$ & $\delta= 0.05$ & $\delta= 0.025$ & $\delta= 0.0125$\\\hline
$\|E_{x}-E_{x}^{\mathrm{ref}}\|_{L^{\infty}}$ & 8.0E-12 & 1.0E-12 & 0 & 0\\
$\|E_{y}-E_{y}^{\mathrm{ref}}\|_{L^{\infty}}$ & 2.0E-12 & 1.0E-12 & 0 & 0\\
$\|E_{z}-E_{z}^{\mathrm{ref}}\|_{L^{\infty}}$ & 1.0E-12 & 0 & 0 &
0\\\hline\hline
\end{tabular}
\end{table}

From Table \ref{table:sol1} it can be seen that without the correction terms,
the solution of VIE has an obvious dependence on the choice of $\delta$ and
the differences follow the order of $O(\delta^{2})$, while the solution is
indeed $\delta$-independent when the correction terms are included, as shown
in Table \ref{table:sol2}.

\subsection{$p$-convergence of the VIE}

For numerical convergence, one can discretize the computational domain into
finer elements (increasing the number $N,$ $h$-refinement) or use higher order
polynomial basis (increasing the number $M,$ $p$-refinement) for the integral
in the VIE. In the current work we focus on the latter since the emphasis is
on the accurate calculation of the Cauchy principal value of the integral with
singularities in the domain.

In the following tests, the Lagrange interpolation of the solution of the VIE
with $M=7$ (namely, 6th order polynomial basis functions) is taken as the
reference solution $\mathbf{E}^{\mathrm{ref}}$ and differences between the
solutions $\mathbf{E}^{M},M=3,4,5,6$ are taken. We consider the linear
relation between the $\log_{10}$ of the energy error
\begin{equation}
\text{Error}=\Vert\mathbf{E}^{M}-\mathbf{E}^{\mathrm{ref}}\Vert_{L^{2}%
(\Omega)},\nonumber
\end{equation}
and the order $p=M-1$ of polynomial basis function,  i.e.
\begin{equation}
\log_{10}(\text{Error})=kp+b,\nonumber
\end{equation}
where $k$ and $b$ are parameters. \begin{figure}[ptb]
\begin{center}%
\begin{tabular}
[c]{ccc}%
\includegraphics[width=0.33\textwidth]{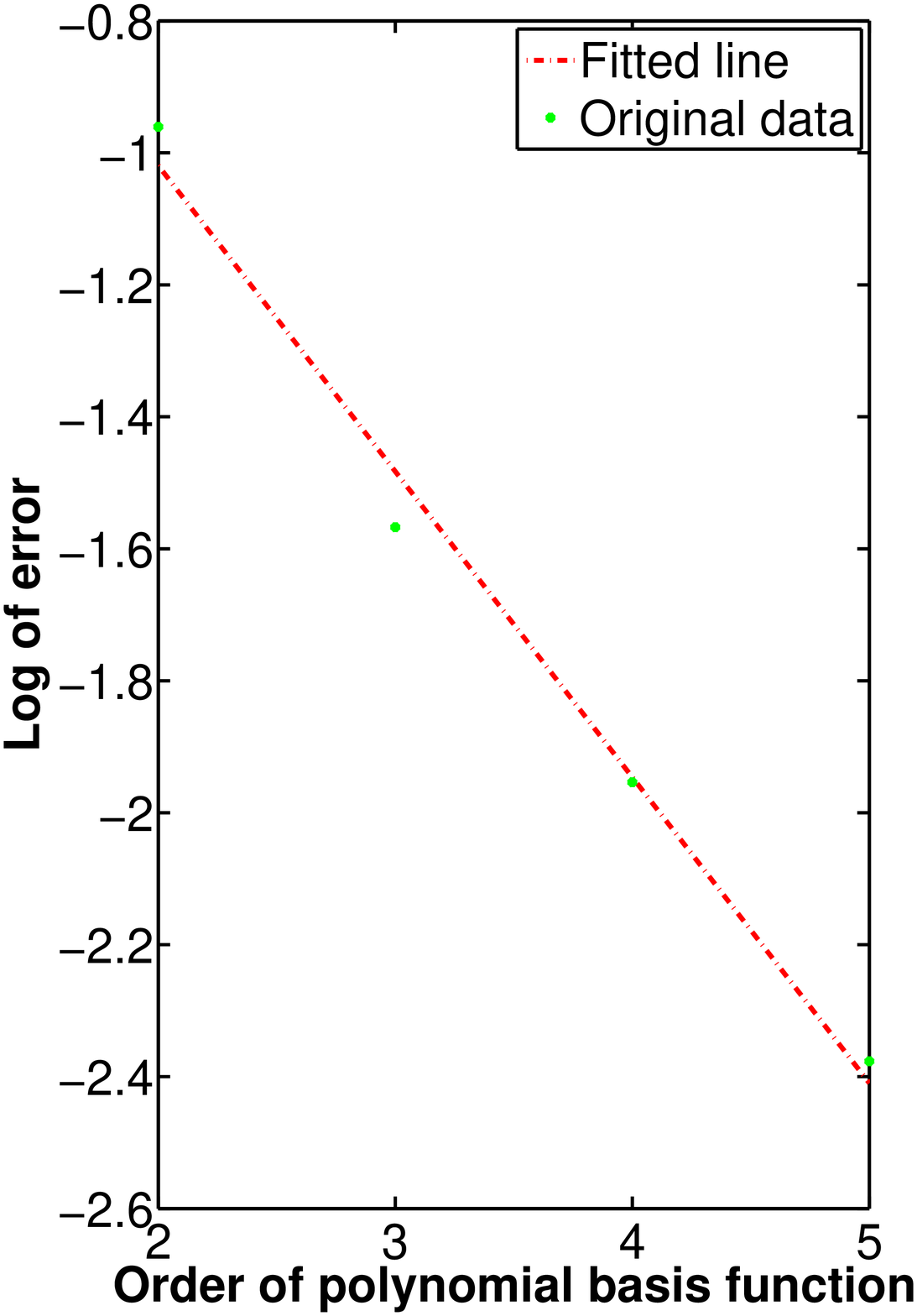} &
\includegraphics[width=0.33\textwidth]{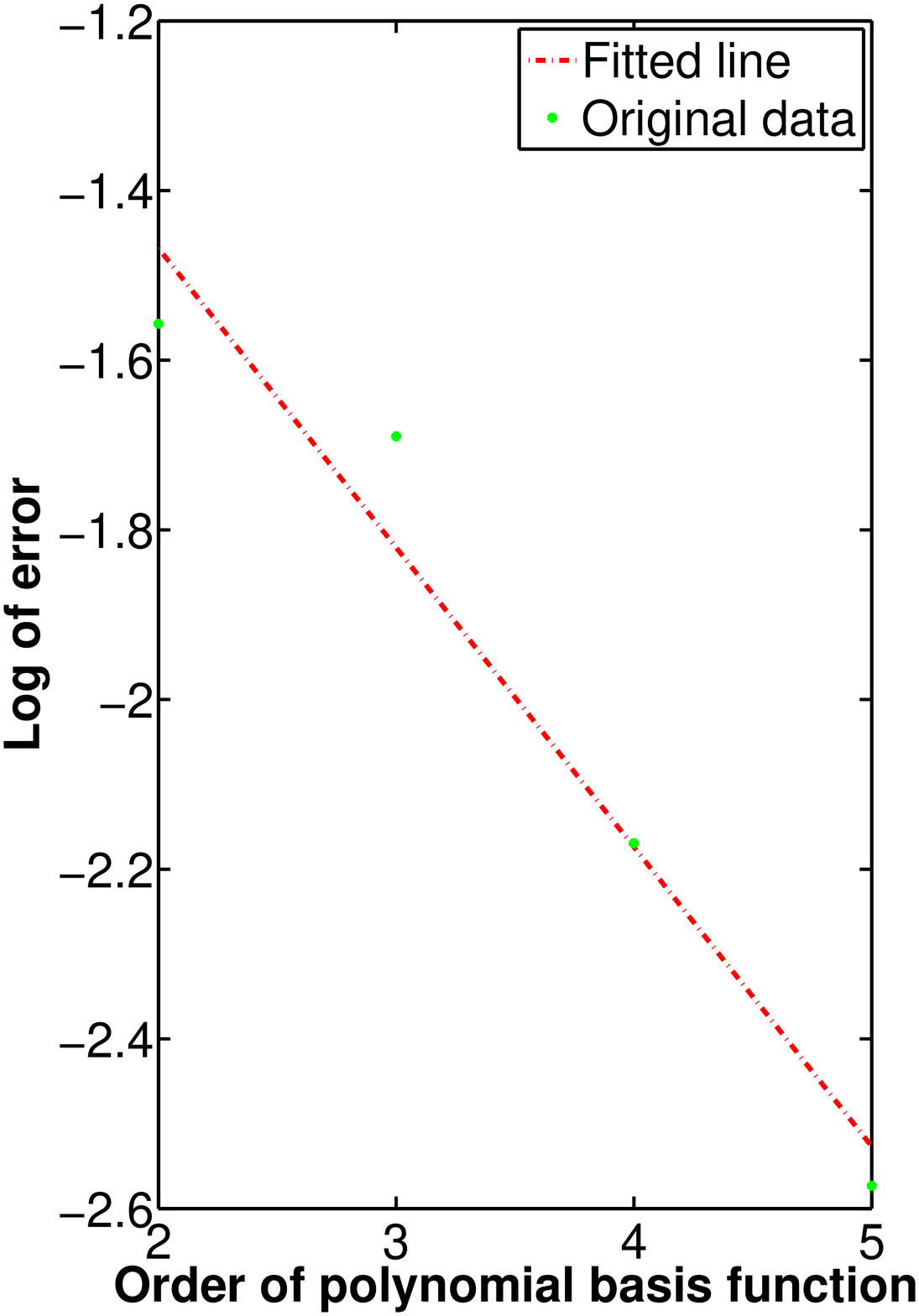} &
\includegraphics[width=0.33\textwidth]{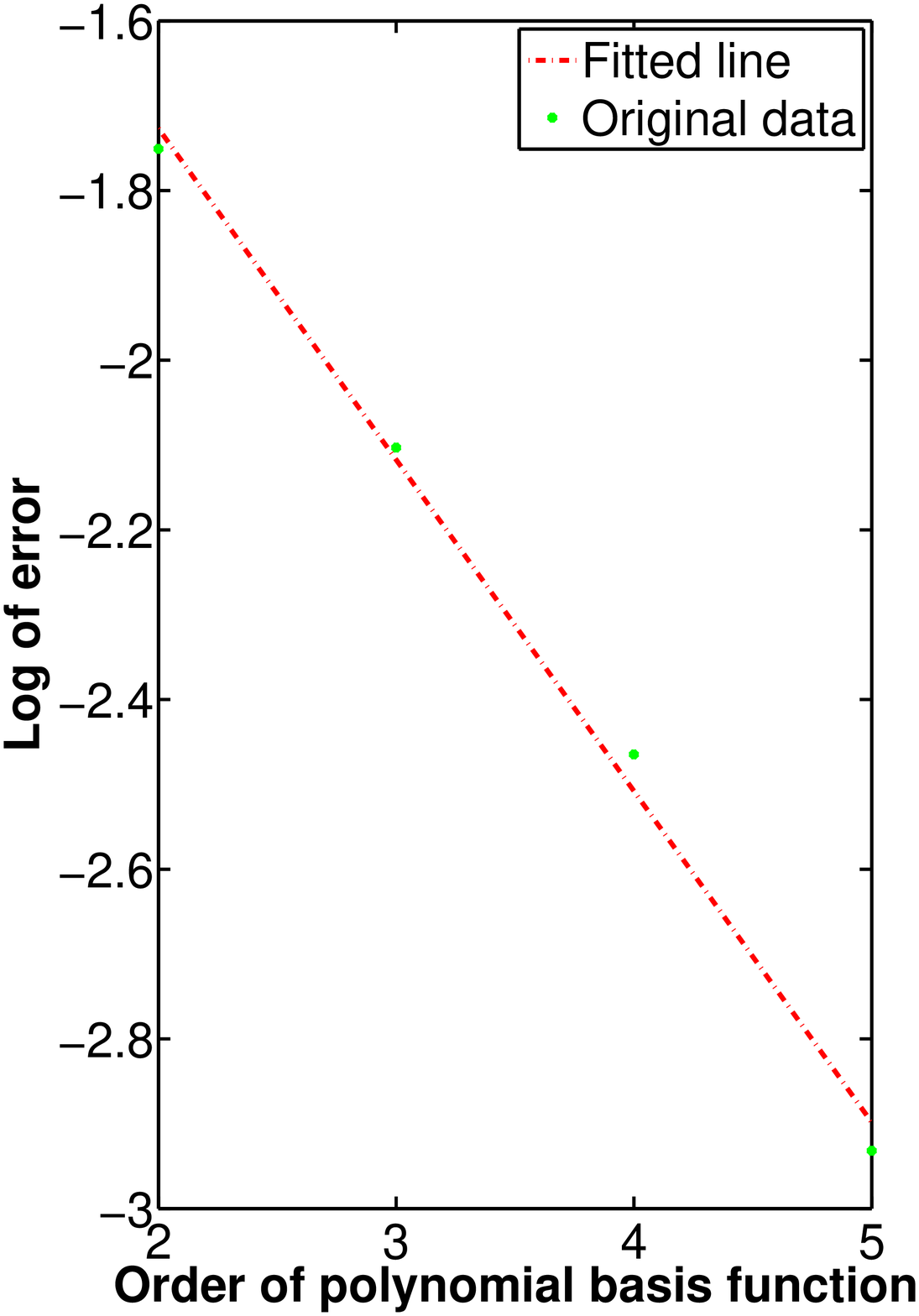}\\
(a) & (b) & (c)\\
&  &
\end{tabular}
\end{center}
\caption{$\log_{10}$ errors of solutions agains order of polynomial basis funtions for $E_{x}$ (a), $E_{y}$ (b) , and $E_{z}$ (c).}%
\label{fig:p-conv}%
\end{figure}
In Fig. \ref{fig:p-conv} there plots the $\log_{10}$ error against
the order $p$, as well as the fitted lines, for the solution components
$E_{x}$, $E_{y}$, and $E_{z}$ respectively. It can be seen that the $\log
_{10}(\mathrm{Error)}$ decays linearly with respect to $p$.
Quantitatively, we obtain
\begin{equation}%
\begin{array}
[c]{ll}%
k=-0.464,\quad b=-0.092, & \quad\text{for}\quad E_{x}\\
k=-0.353,\quad b=-0.763, & \quad\text{for}\quad E_{y}\\
k=-0.391,\quad b=-0.946, & \quad\text{for}\quad E_{z}.
\end{array}
.\nonumber
\end{equation}
Hence we conclude that the exponential convergence is achieved for the errors
against the order $p$.

\begin{figure}[ptb]
\begin{center}%
\begin{tabular}
[c]{cccc}%
\includegraphics[width=0.23\textwidth]{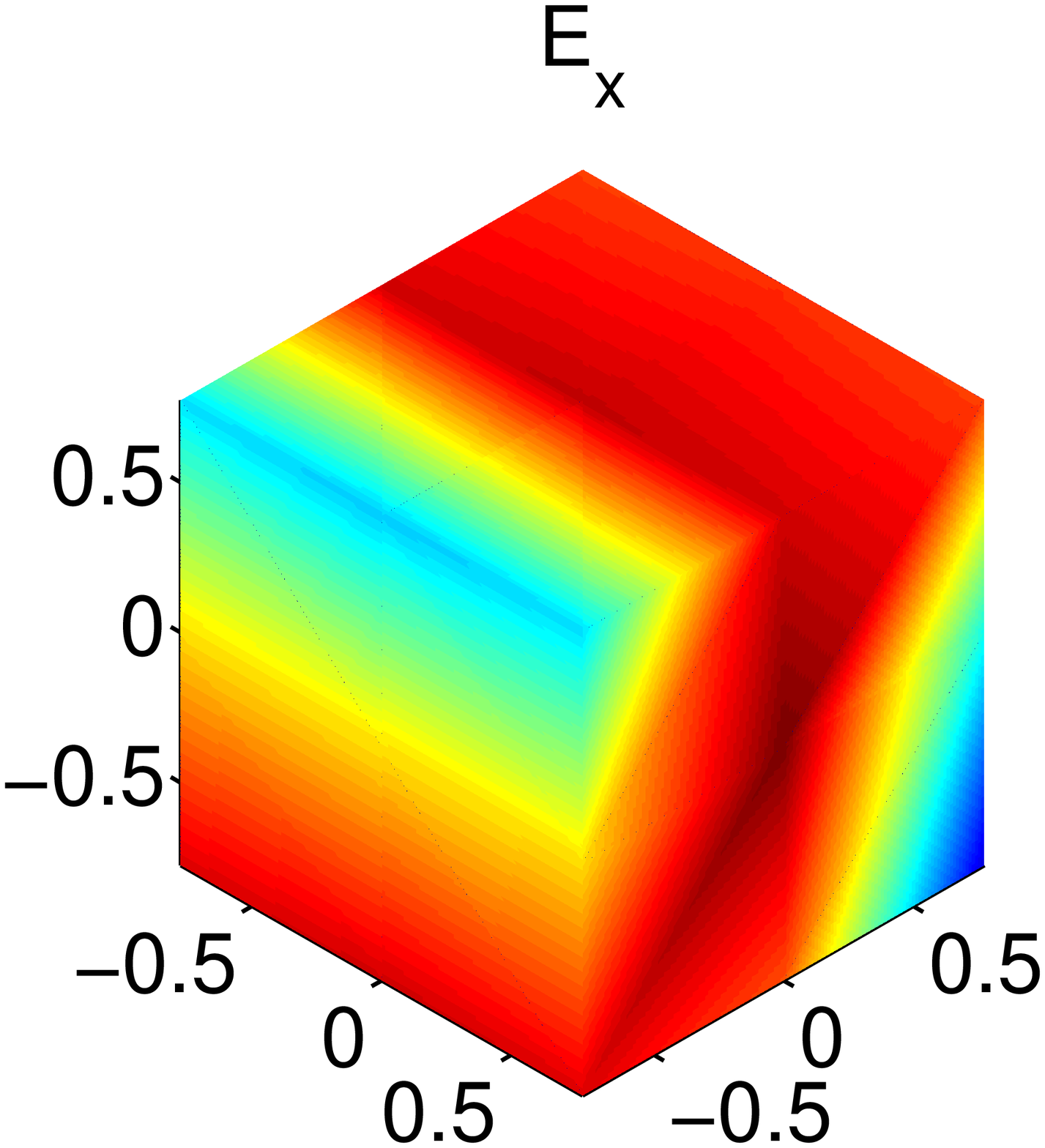} &
\includegraphics[width=0.23\textwidth]{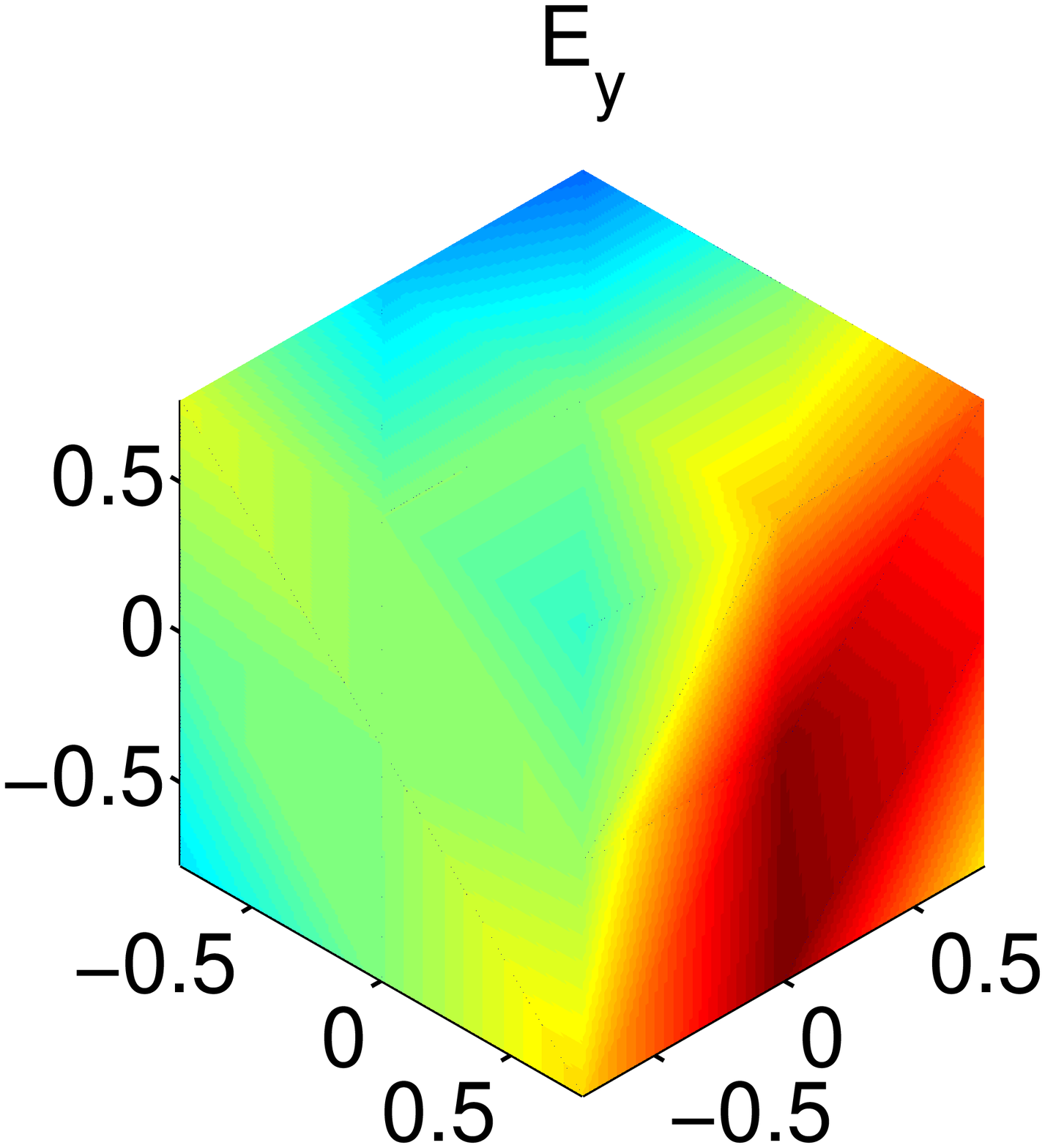} &
\includegraphics[width=0.23\textwidth]{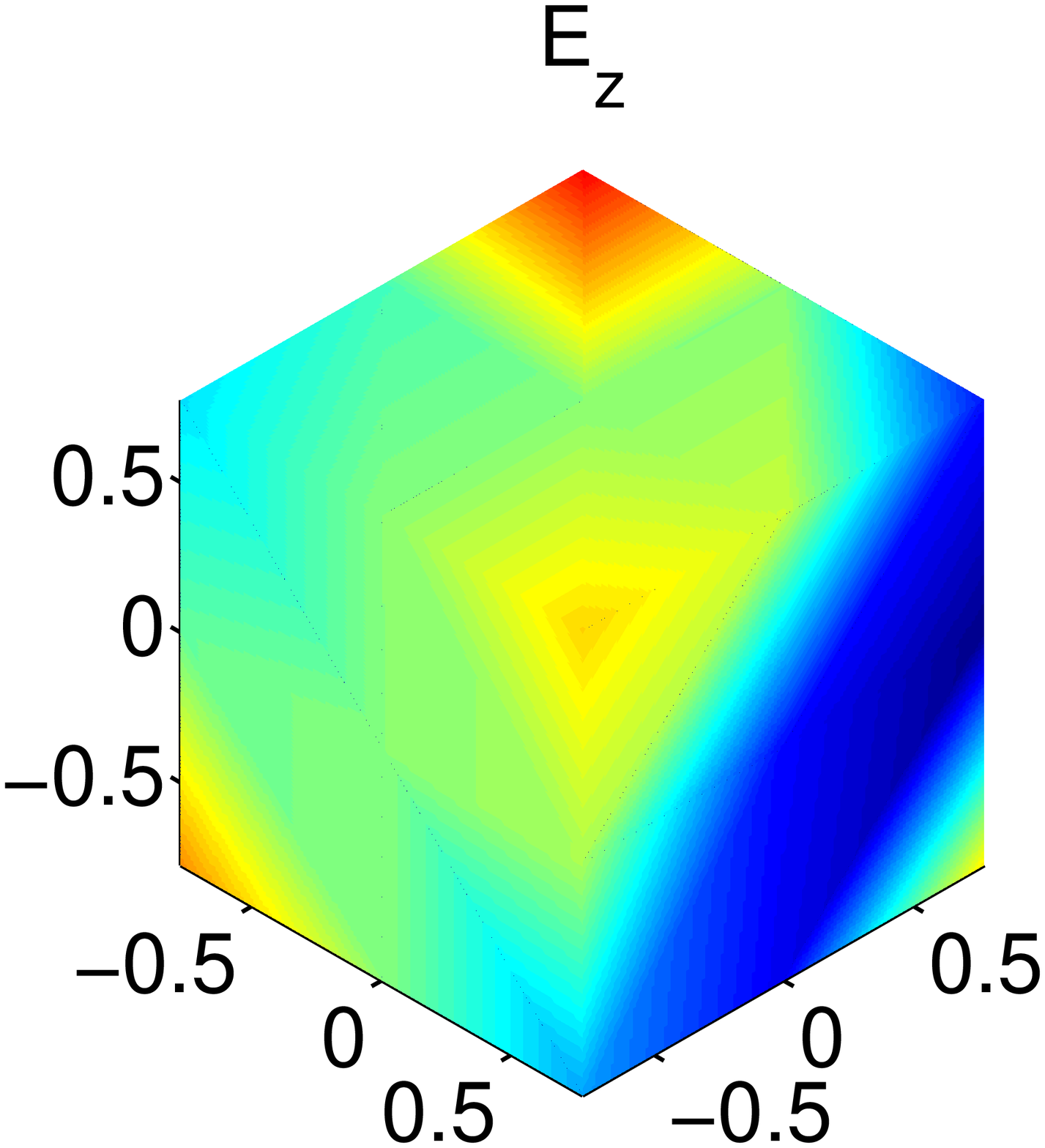} &
\includegraphics[width=0.23\textwidth]{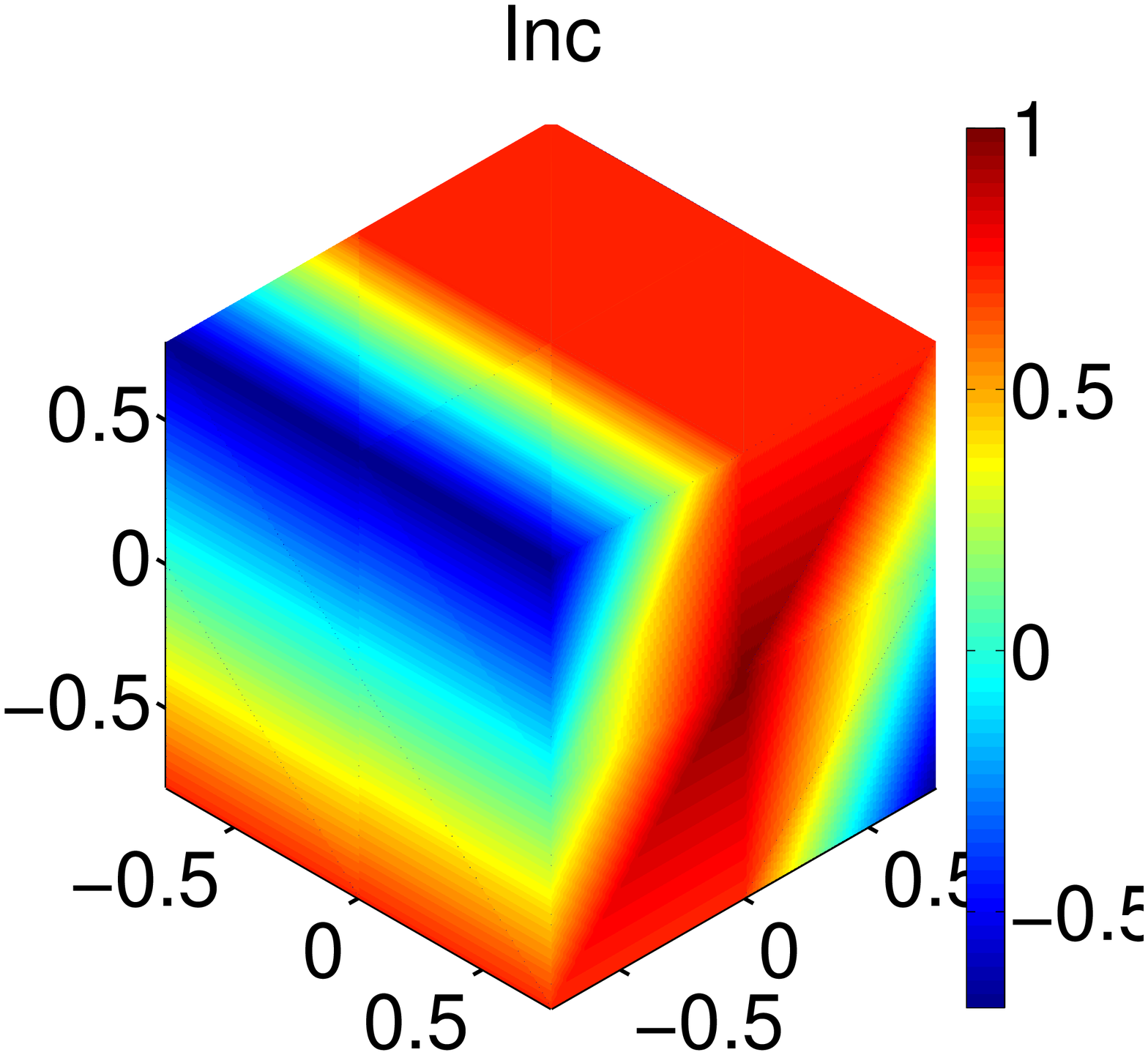}\\
(a1) & (b1) & (c1) & (d1)\\
\includegraphics[width=0.23\textwidth]{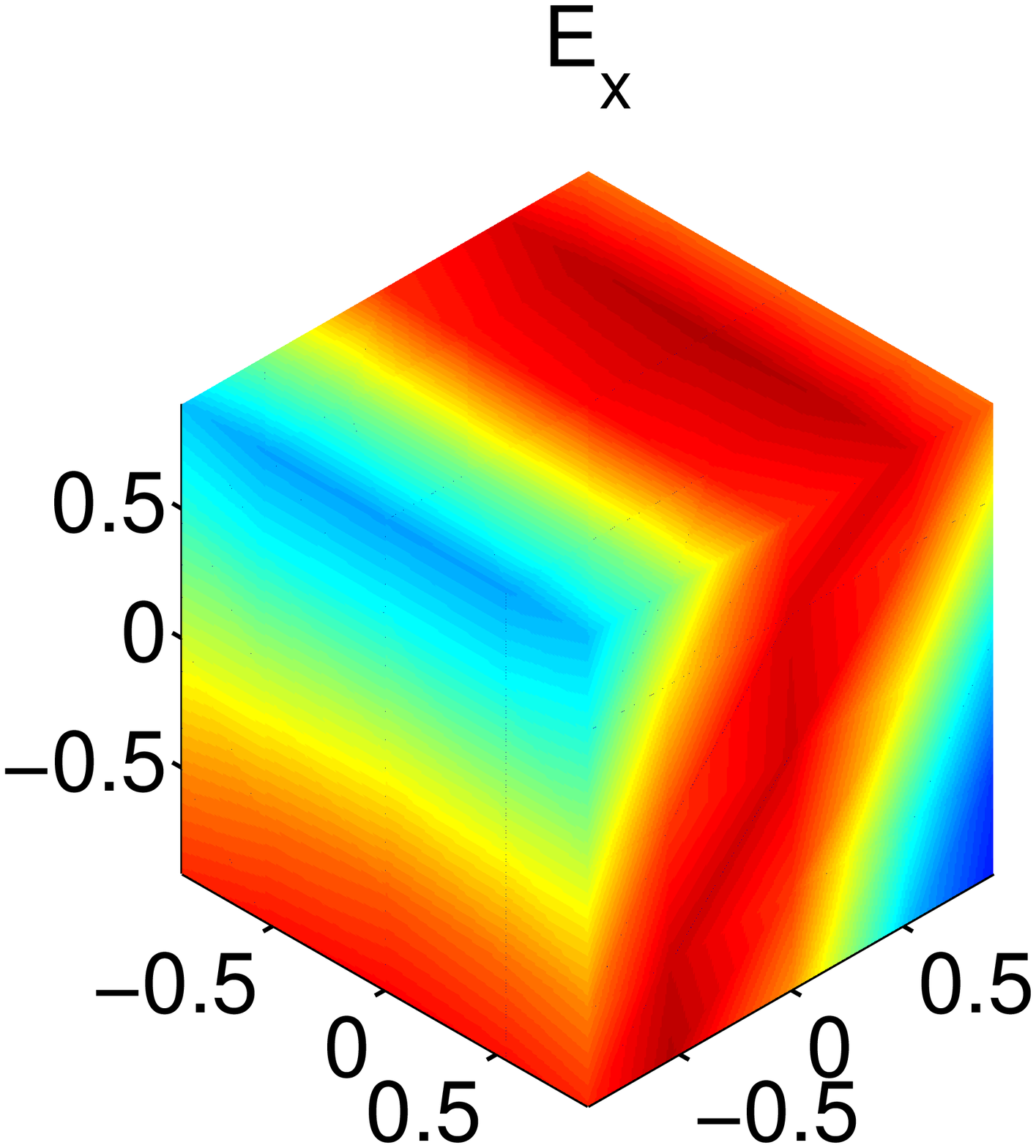} &
\includegraphics[width=0.23\textwidth]{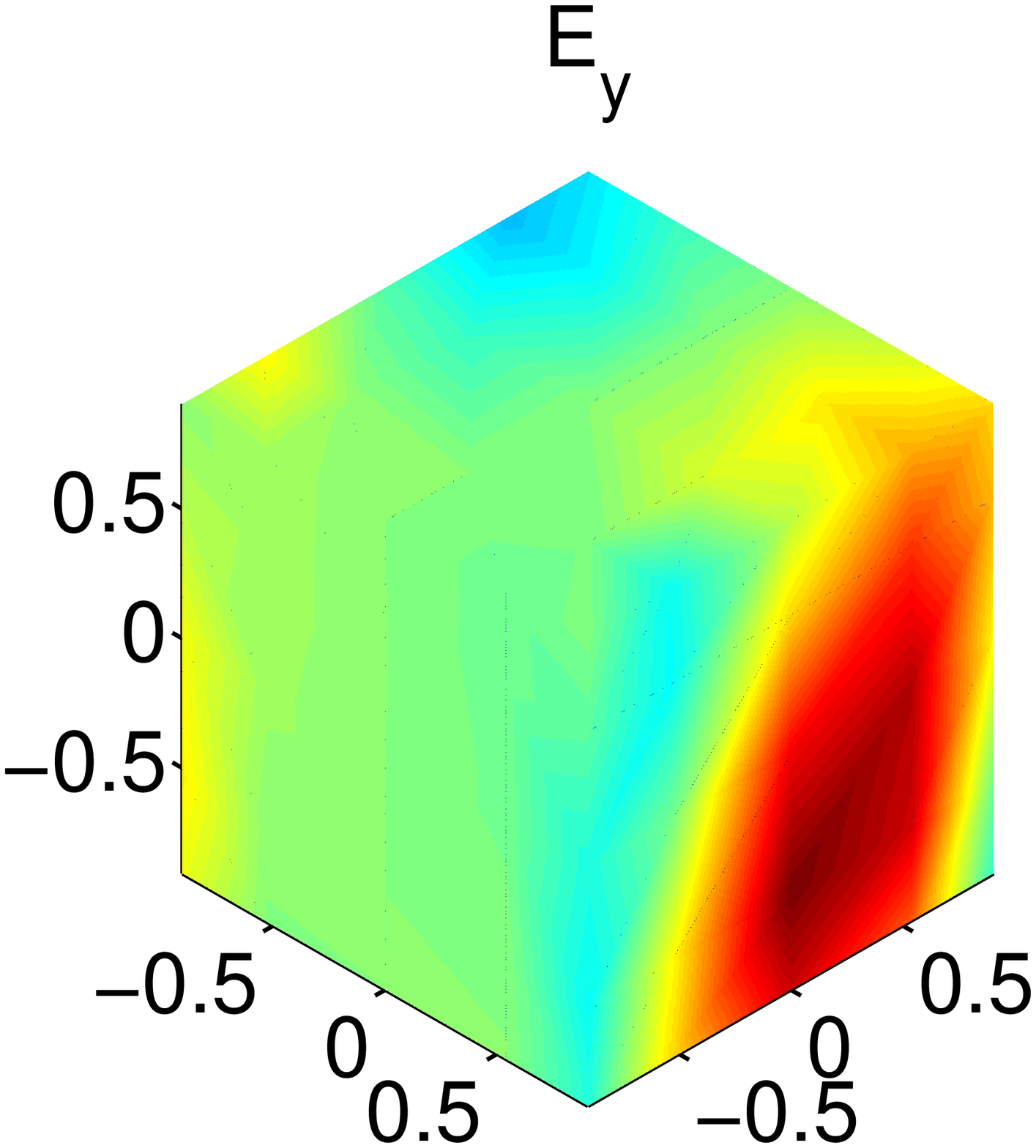} &
\includegraphics[width=0.23\textwidth]{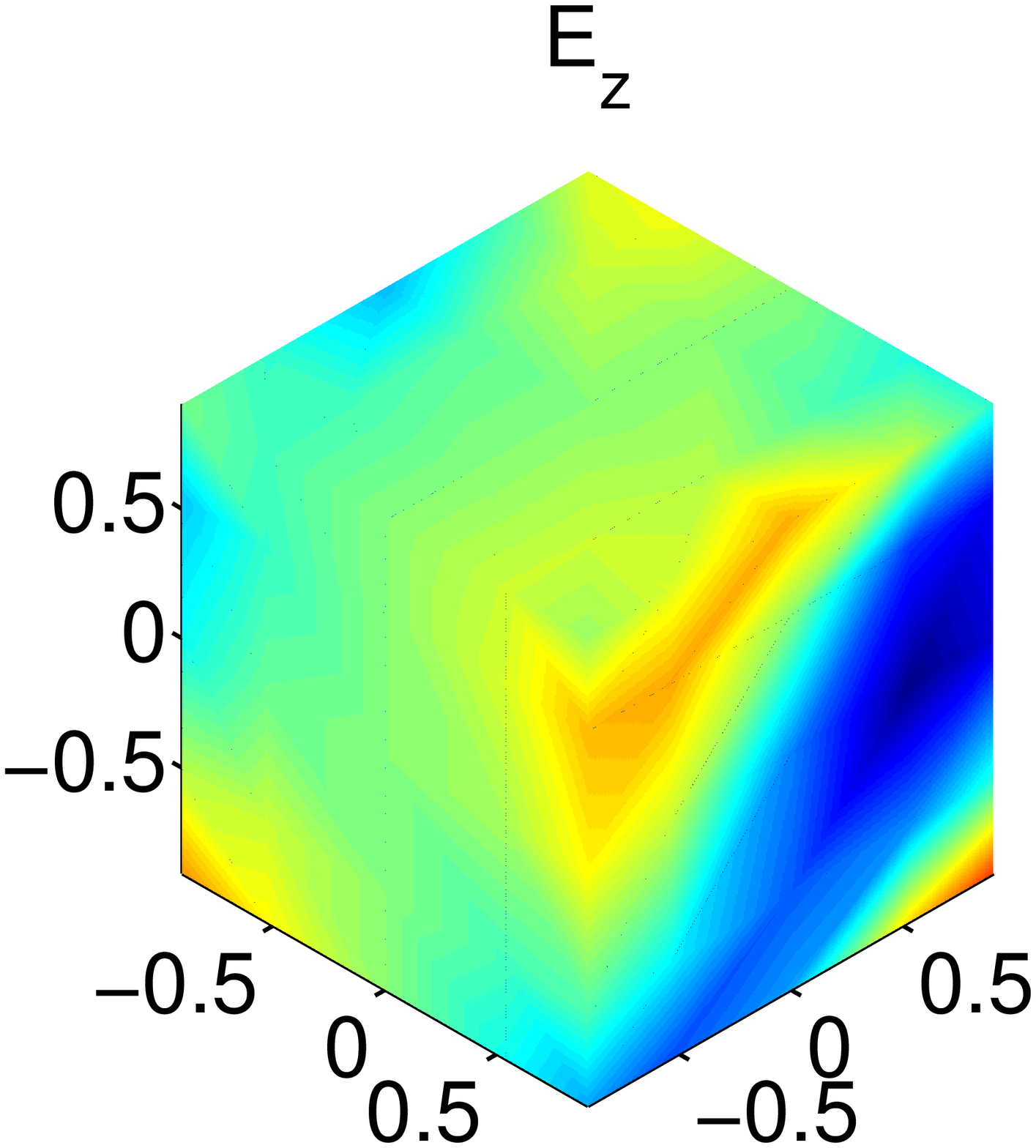} &
\includegraphics[width=0.23\textwidth]{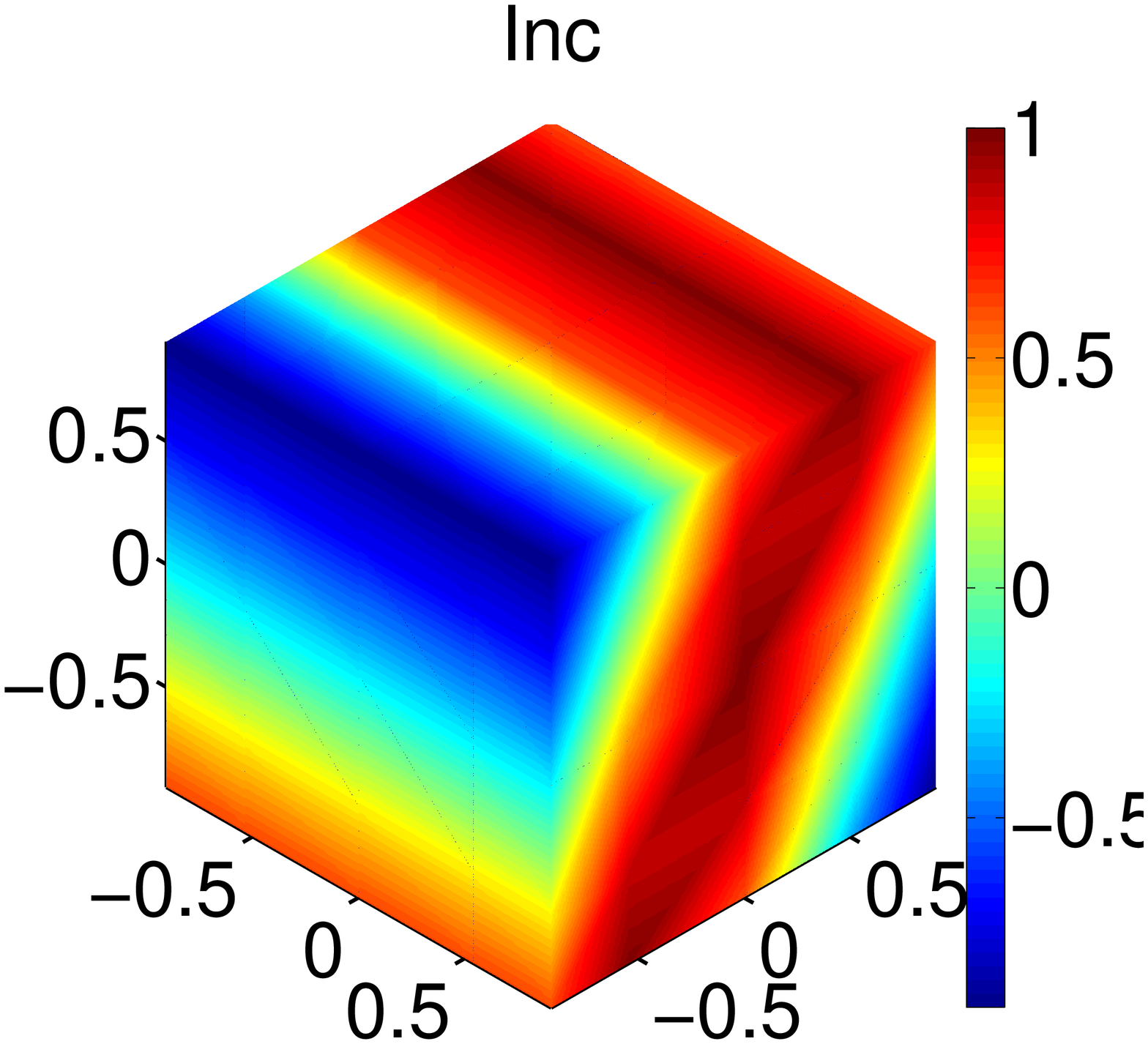}\\
(a2) & (b2) & (c2) & (d2)\\
\includegraphics[width=0.23\textwidth]{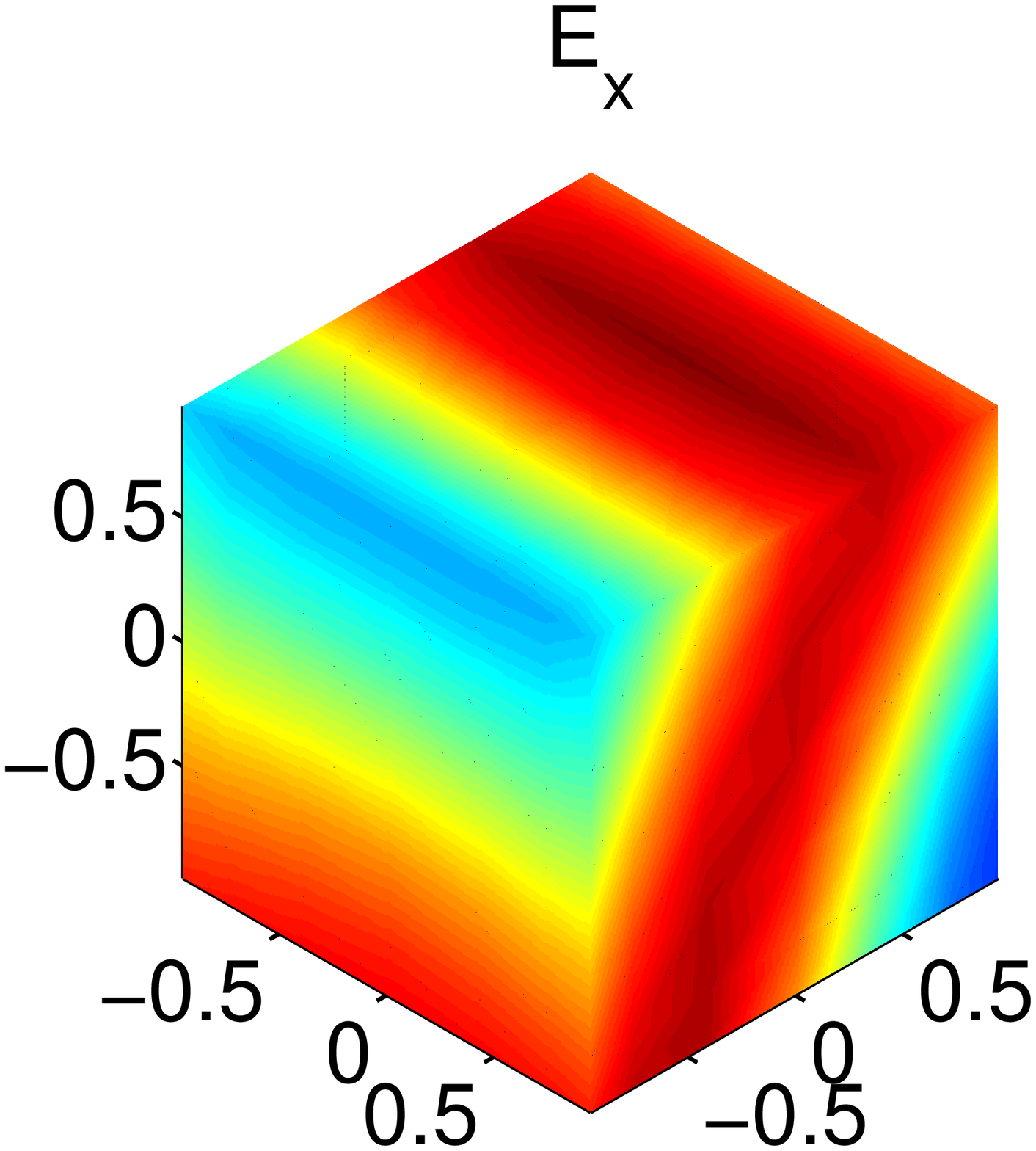} &
\includegraphics[width=0.23\textwidth]{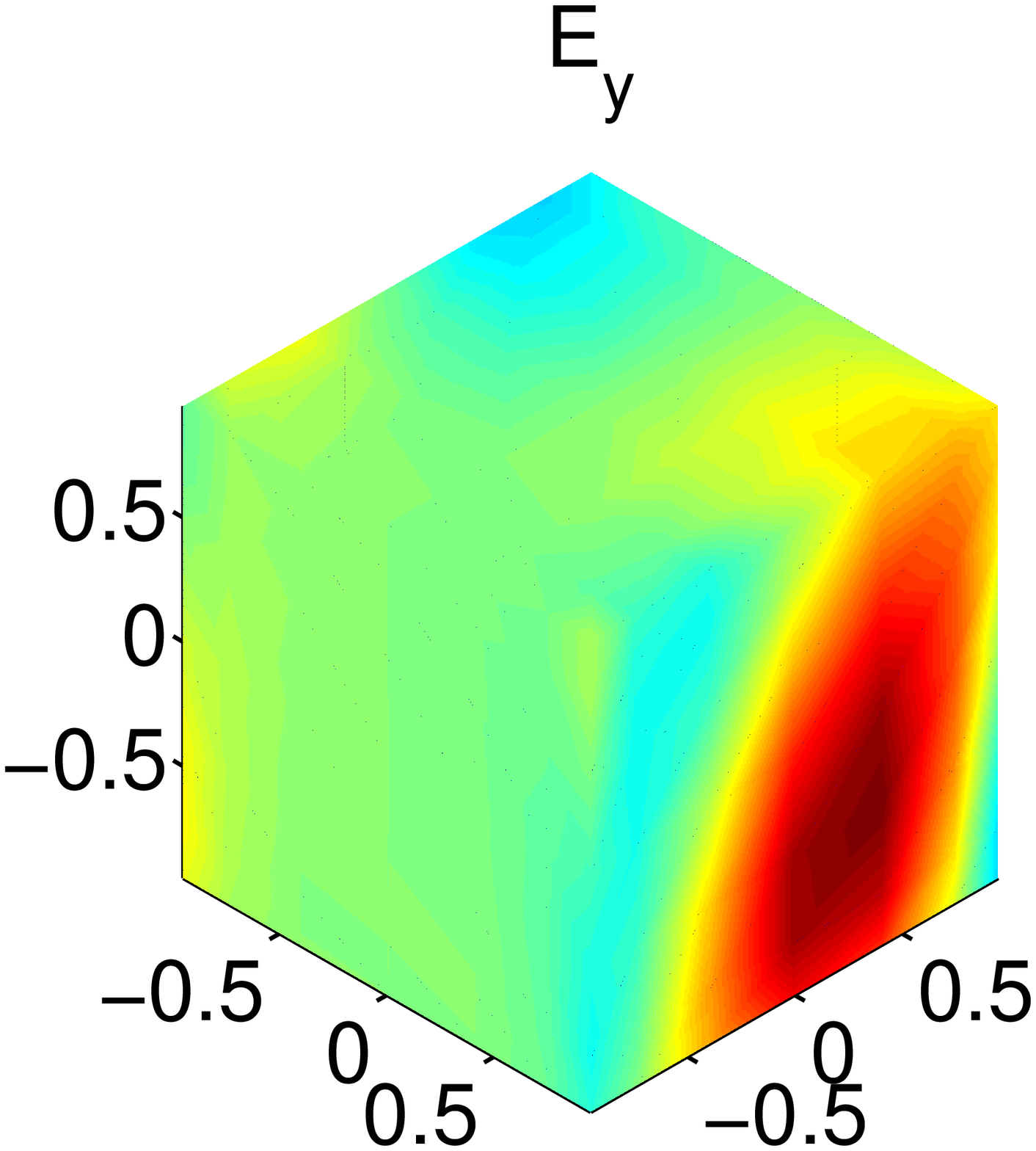} &
\includegraphics[width=0.23\textwidth]{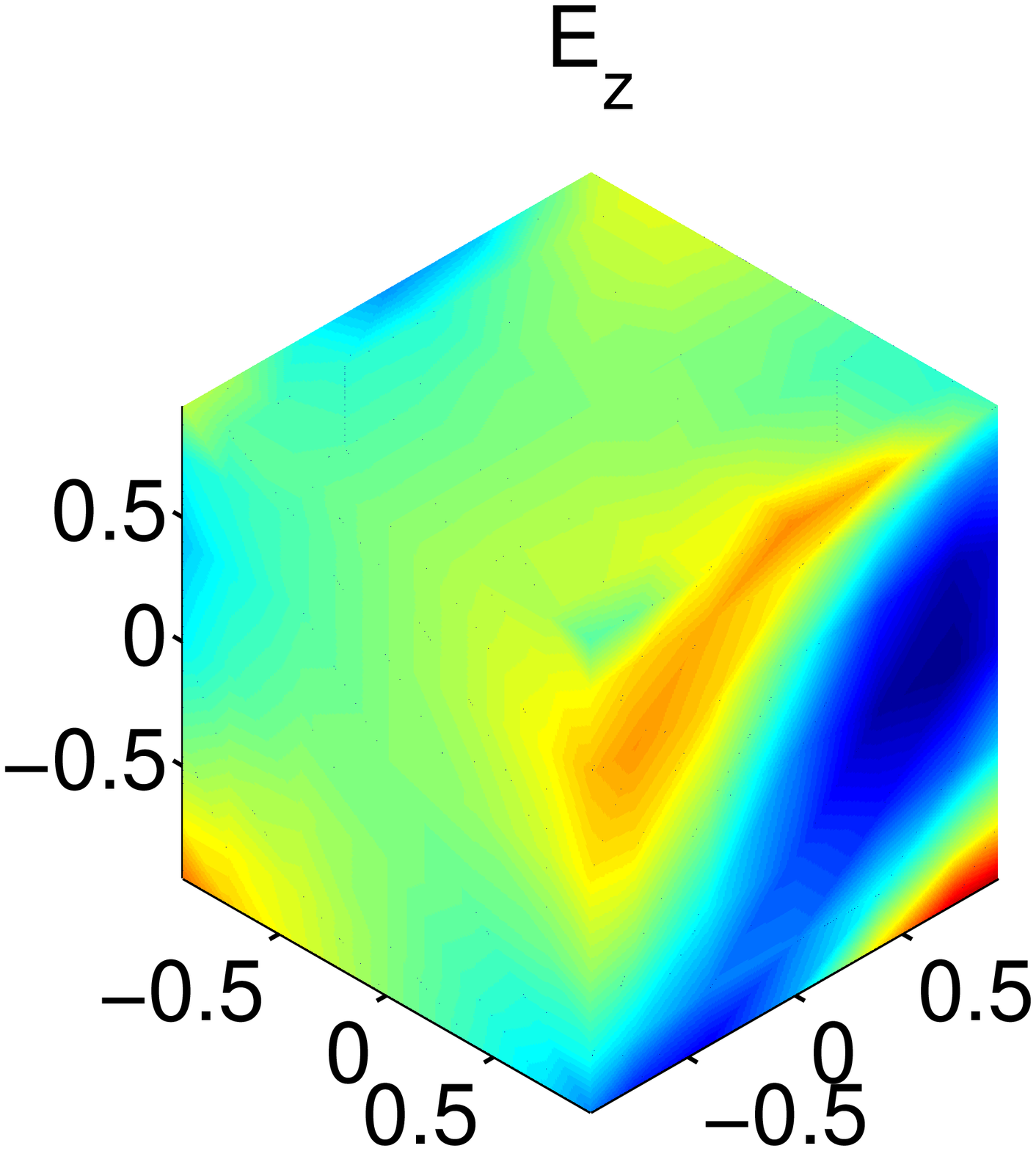} &
\includegraphics[width=0.23\textwidth]{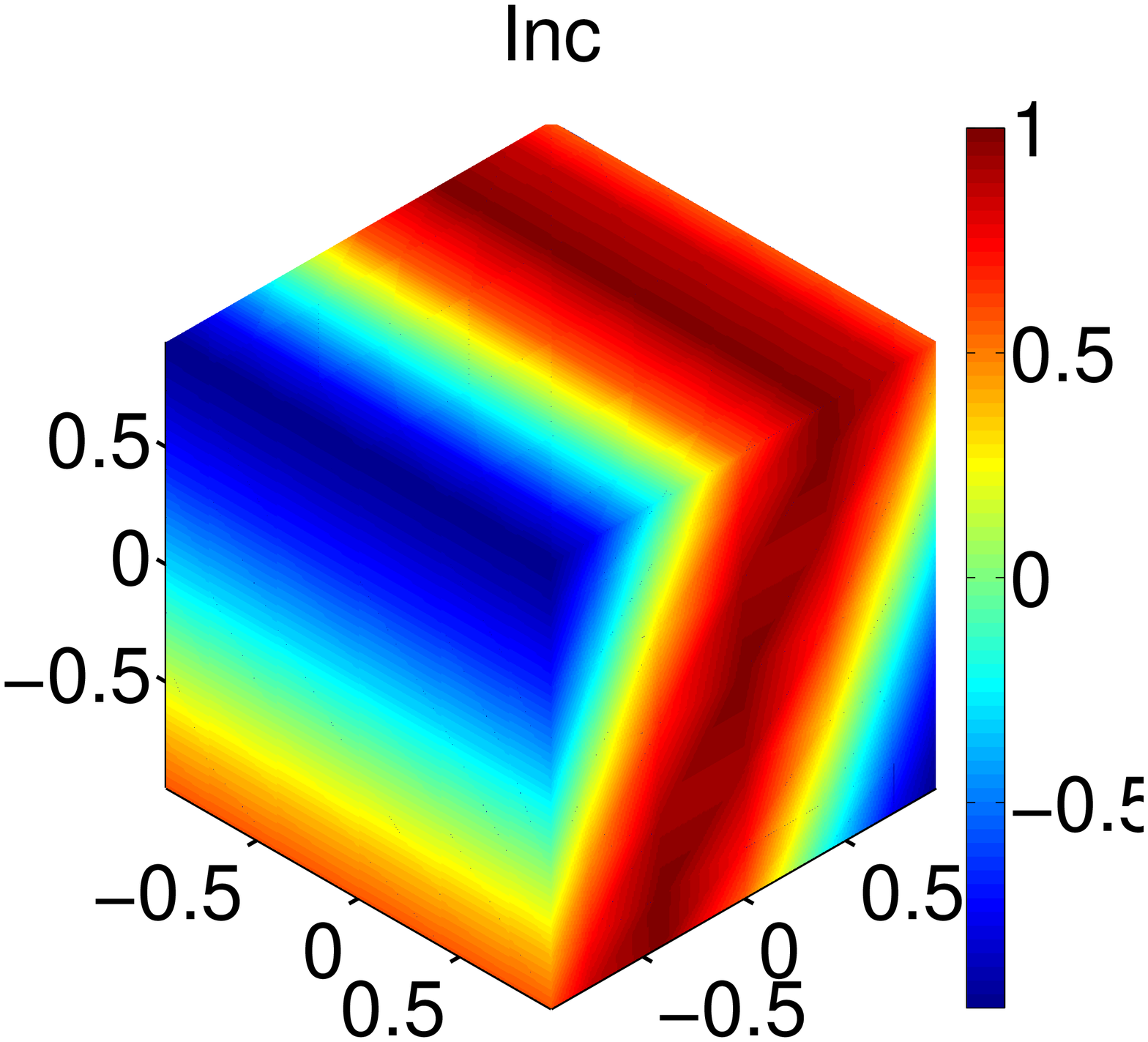}\\
(a3) & (b3) & (c3) & (d3)\\
&  &  &
\end{tabular}
\end{center}
\caption{Electric field and the incident wave with $M=3$ (first row), $M=5$
(second row), and $M=7$ (third row).}%
\label{fig:p-refine}%
\end{figure}

Figure \ref{fig:p-refine} is a 3D plot of the incident wave and resulting
electric fields inside the cubic scatter. Rows 1-3 of Fig. \ref{fig:p-refine}
are for the computations with Gauss node $M=3, 5, 7$, respectively, and from
the left to the right are the electric field components $E_{x}$, $E_{y}$,
$E_{z}$ and the incident wave $\mathbf{E}_{x}^{\mathrm{inc}}$, respectively.

\begin{figure}[ptb]
\begin{center}%
\begin{tabular}
[c]{cc}%
\includegraphics[width=0.5\textwidth]{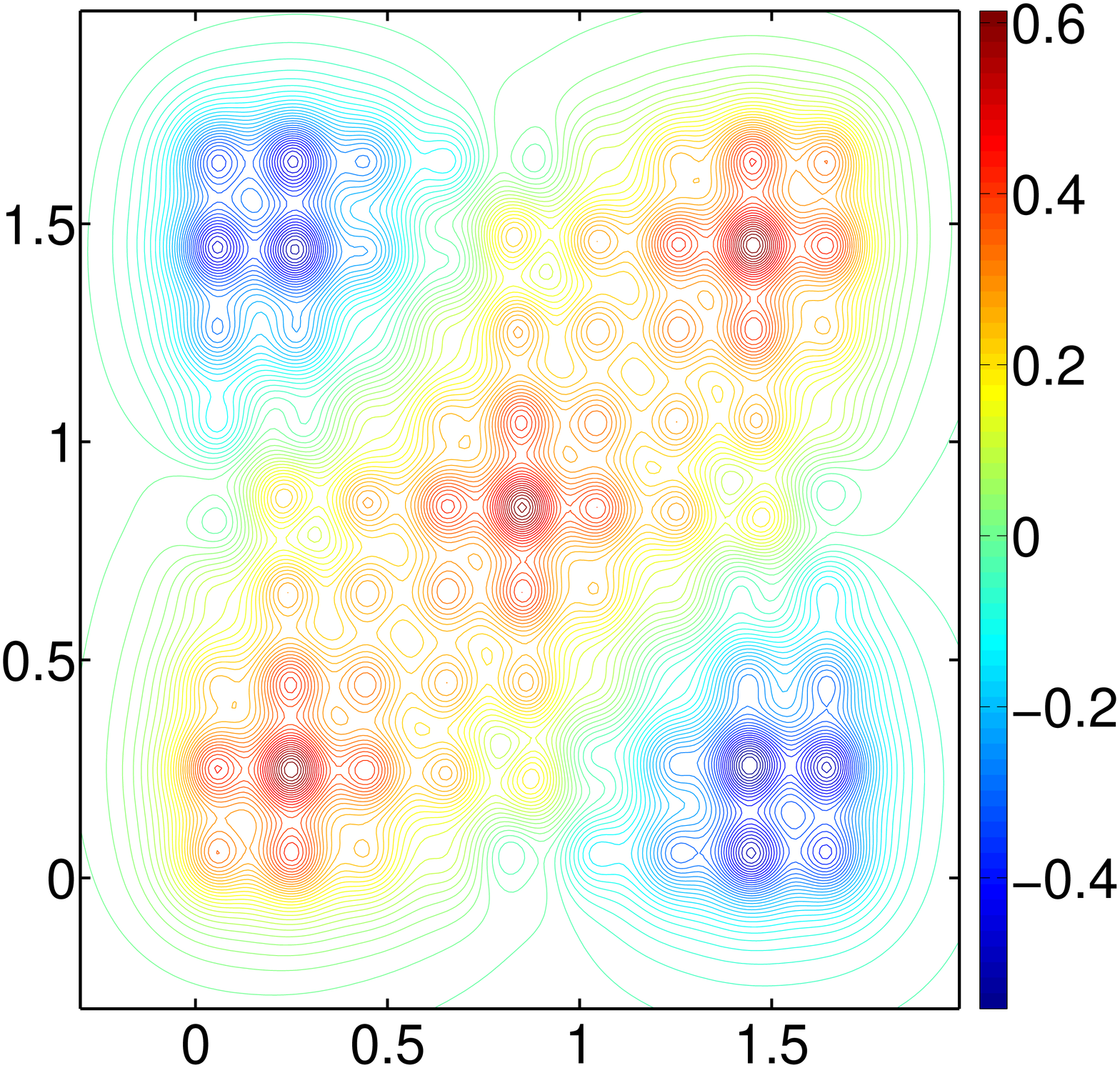} &
\includegraphics[width=0.5\textwidth]{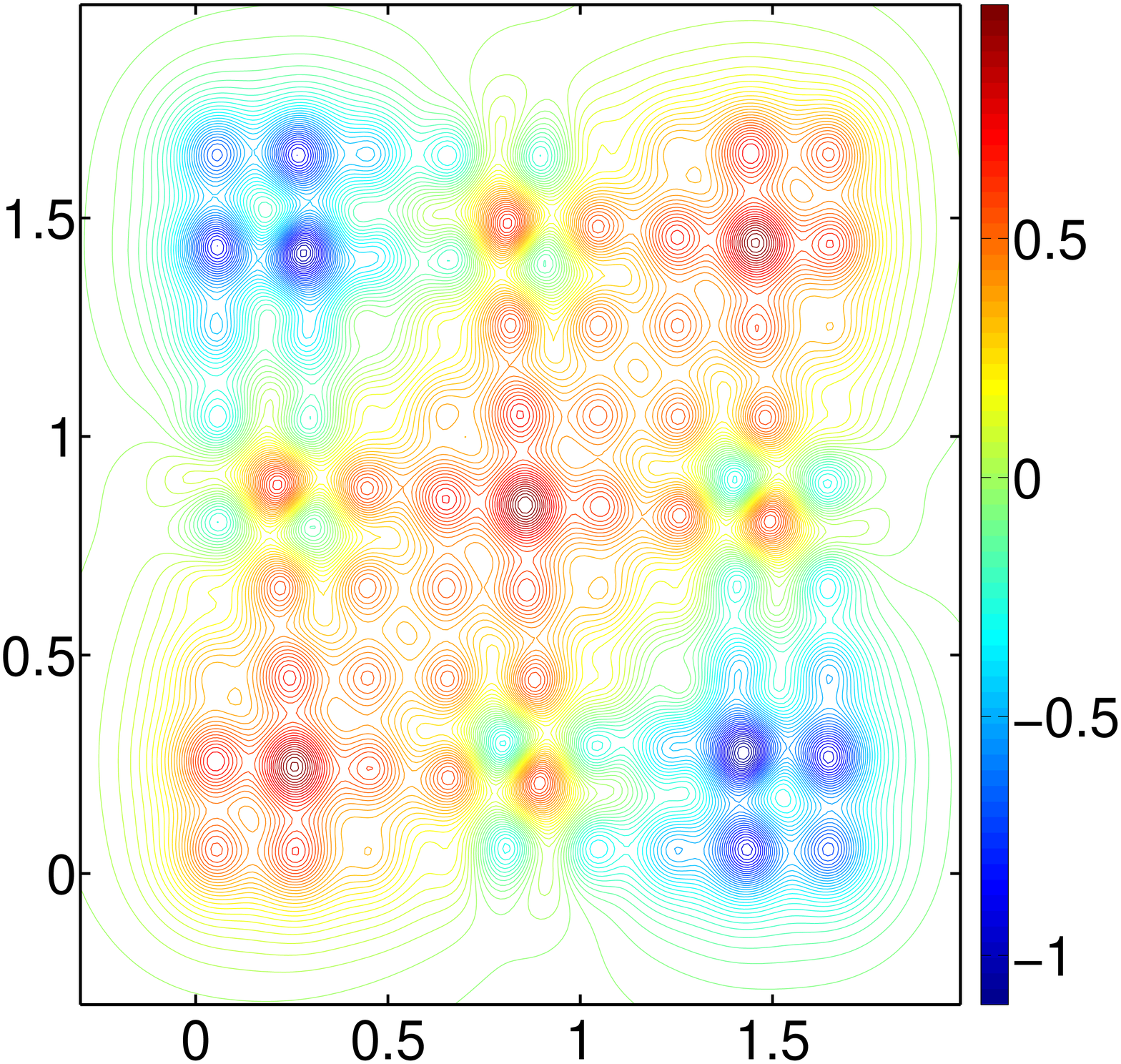}\\
(a) & (b)\\
\includegraphics[width=0.5\textwidth]{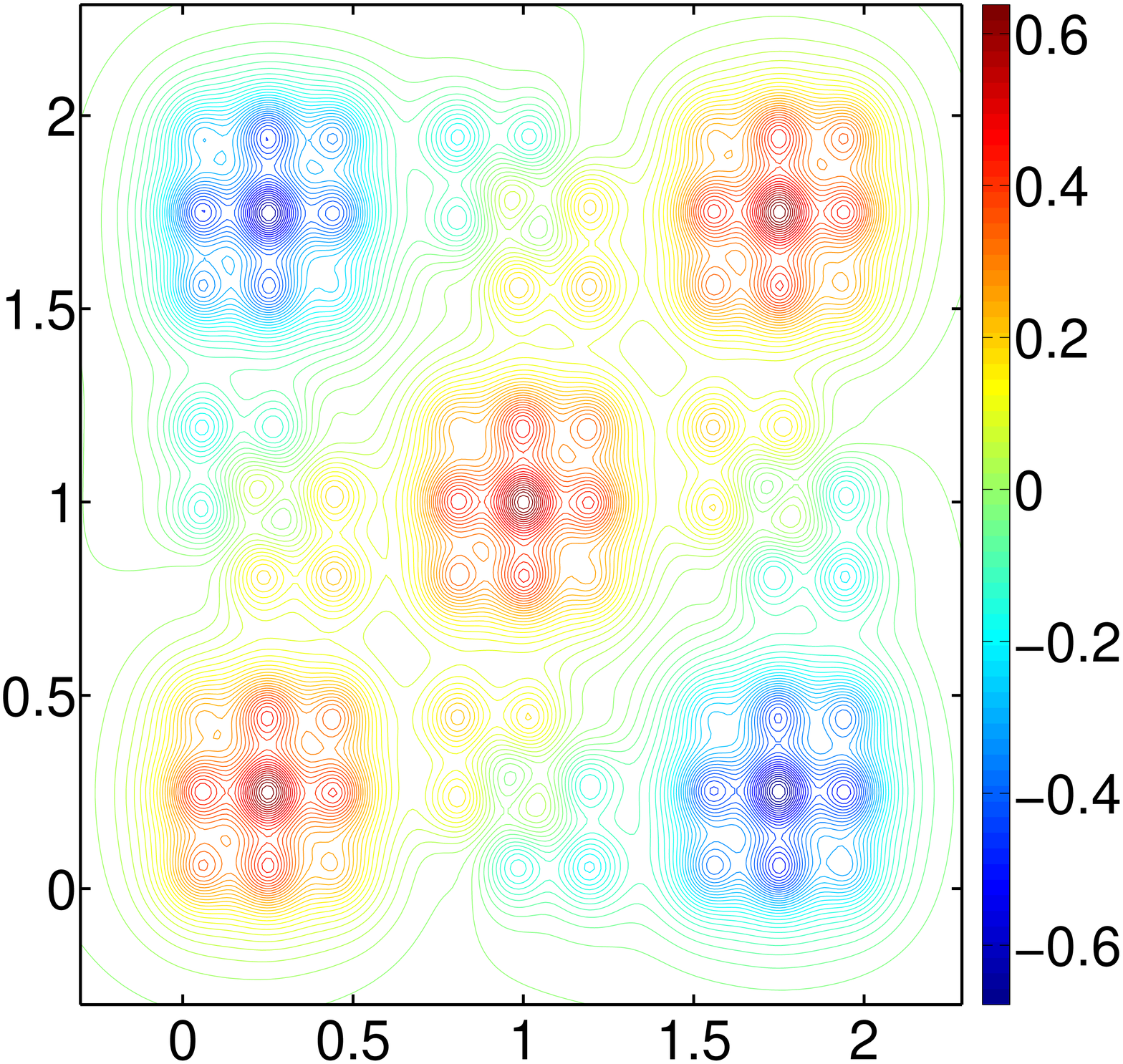} &
\includegraphics[width=0.5\textwidth]{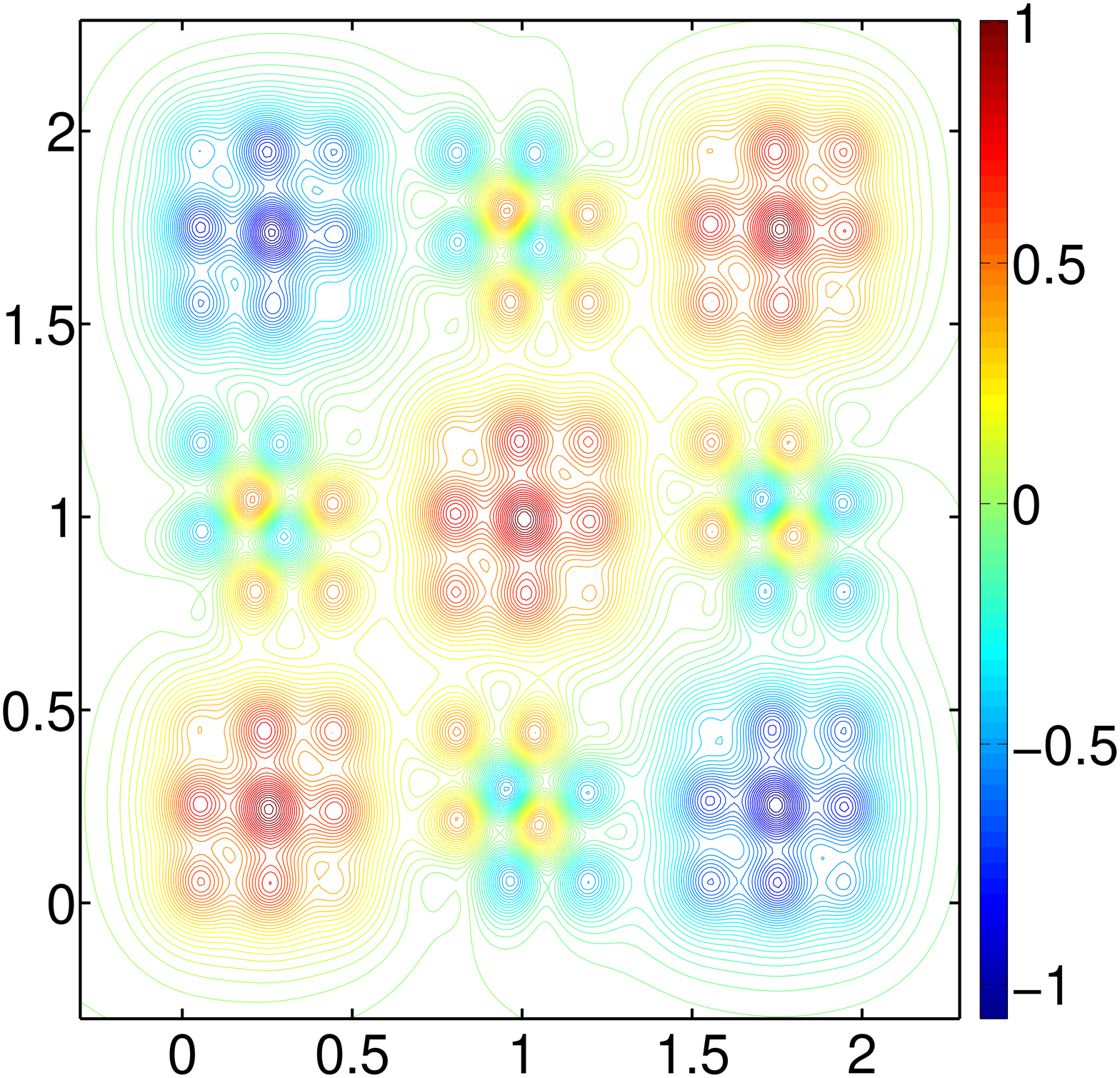}\\
(c) & (d)
\end{tabular}
\end{center}
\caption{Contours of the electric field from scatter arrays. (a)-(b): scatters
are 0.1 apart from each other; (c)-(d): scatters are 0.25 apart from each
other. (a, c): $\Delta\epsilon=4$; (b, d): $\Delta\epsilon=16$. }%
\label{fig:array}%
\end{figure}

Figure \ref{fig:array} displays the electric field in the free space where
nine cubic scatters are embedded. In order to clearly illustrate the wave
interactions among the scatter array, only the scattering field $\mathbf{E}%
(\mathbf{r)-\mathbf{E}^{\mathrm{inc}}(r)}$ are shown in contours in a specific
cross-section plane. In these tests, the incident wave is taken as
\begin{equation}
\mathbf{E}_{x}^{\mathrm{inc}}=\mathbf{E}_{y}^{\mathrm{inc}}=0,\quad
\mathbf{E}_{z}^{\mathrm{inc}}=e^{ik(-2x+2y)}.\nonumber
\end{equation}
Each cube has a length of 0.5 and they form a $3\times3$ array align in the
$x$-$y$ plane. The coordinate of the center of the first cube is
$(0.25,0.25,0.25)$. We consider two cases: in the first case the cubes are
relatively far way from each other (0.25 apart) and in the second case the
cubes are relatively closer (0.1 apart). In each case a lower permittivity
($\Delta\epsilon=4$) and a higher permittivity ($\Delta\epsilon=16$) are
taken. Contours of the electric field ${E}_{x}$ in the plane $z=-0.1$ are
shown in Fig {\ref{fig:array}}. The first and second rows are for the scatters with smaller and bigger distances, respectively, in which the left is for
$\Delta\epsilon=4$ and the right is for $\Delta\epsilon=16$. It can be seen
that the electric waves have stronger interactions when scatters are closer to
each other.

\section{Conclusion}

\label{sec:conclusion}

In this paper, we have developed an highly accurate and efficient method to
solve numerically the volume integral equation (VIE) of the Maxwell equation
involving the Cauchy principal value of the singular dyadic Green's function.
The VIE used allows us to compute the CPV for finite size exclusion volume
$V_{\delta}$ without truncation errors (of the order of $O(\delta)$ in general
and it is $O(\delta^{2})$ when the excluded volume is a sphere), which were
accurately accounted for by including corrections terms. With the correction
terms, the numerical solution of the VIE will be $\delta$-independent. In
addition, an efficient quadrature formula has been introduced to accurately
compute the ${\frac{1}{R^{3}}}$ type integration on the domain $\Omega
\backslash V_{\delta}$. Finally, a Nystr\"{o}m collocation method based on the
proposed quadrature formula were applied to the VIE resulting in an accurate
and $\delta$-independent solution of the electric field.

The developed algorithms are verified with several numerical tests. First, the
accuracy of the interpolated quadrature weights for dyadic Green's function
was confirmed for various locations of the spherical exclusion volume
$V_{\delta}$ inside a cubic scatter. Secondly, the procedure of computing the
Cauchy principal values for the dyadic Green's function with a finite size
$\delta$ is shown to be accurate with the help of the proposed correction
terms. As a result, the solution of the VIE is shown to be $\delta
$-independent, and finally, we demonstrated the convergence of the VIE for
$p$-refinement with increasing order of basis functions $p=m-1$, where the
number of Gauss node along each direction $m=3,4,5,6$, and $7$, respectively.

\bibliographystyle{plain}
\bibliography{vie}

\end{document}